\documentclass[11pt]{article}

\usepackage[utf8]{inputenc}
\usepackage[T1]{fontenc}
\usepackage[margin=1in]{geometry}

\usepackage{amsmath,amssymb,amsfonts}
\usepackage{graphicx}
\usepackage{caption}
\usepackage{xcolor}
\usepackage{url}
\usepackage[hidelinks]{hyperref}

\usepackage[framed,numbered,autolinebreaks,useliterate]{mcode}

\title{Mathematical modelling and simulation of HIV dynamics}
\date{}
\author{Abdul Rab\\
Department of Mathematics and Social Science\\
Sukkur IBA University, Pakistan\\
\texttt{abdulrab@live.com}} 
\begin{document}

\maketitle

\begin{abstract}
We study within-host HIV dynamics using a three--component nonlinear ordinary differential equation model for healthy CD4$^{+}$ T cells, infected CD4$^{+}$ T cells, and free virus. In addition to the baseline model without treatment, we consider two treatment extensions that incorporate antiretroviral therapy: (i) separate efficacy terms for Reverse Transcriptase inhibitors and Protease inhibitors acting on infection and virus production, and (ii) a simplified formulation using a single combined efficacy parameter. For each model we determine equilibrium points and apply linearization to obtain the Jacobian matrix and local stability conditions near equilibrium. We perform numerical simulations using the classical fourth--order Runge--Kutta method to illustrate the evolution of cell populations and viral load under different therapy levels and treatment start times, including continuous treatment scenarios. The simulations highlight how therapy strength and timing shape transient behavior and can lead to sustained viral suppression.
\end{abstract}


\section{Mathematical modelling}\label{chapter1}
In this chapter, dynamical models of the HIV are reviewed to
understand how the HIV virus effects the immune system, as it enters the
human body. The effect of antiretoviral treatment i.e. Reverse Transcriptase Inhibitors (RTIs) and Protease Inhibitors (PIs) added in simple HIV dynamical model and by further simplification effect RTIs and PIs is combined as single parameter in the HIV model. In addition, the
system is linearized around the equilibrium point, which leads to the
first order system of linear ordinary differential equations.
Finally, the HIV dynamical models have been numerically simulated and assessed in Chapters~\ref{chapter2}$-$\ref{chapter3}, and results have been interpreted and analysed~\cite{Moysis2016}.

\subsection{Introduction}
\begin{figure}[ht]
    \centering
    \includegraphics[width=5in,height=2in]{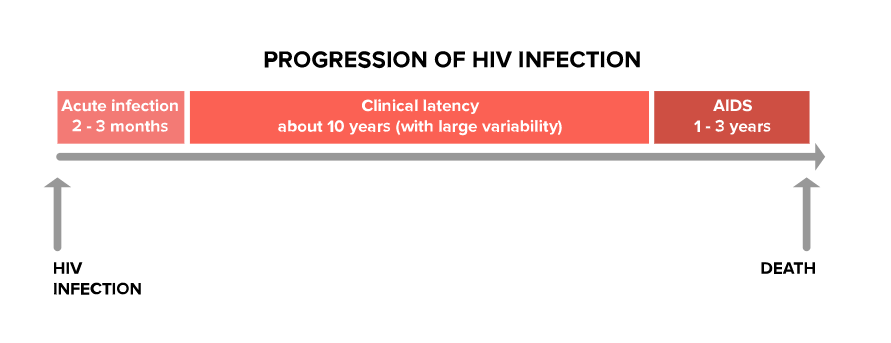}
    \caption{Source: Khan Academy (Course on HIV Infection)~\cite{khan_academy}.}
    \label{fig:figure_1}
\end{figure}

\textbf{What is HIV?} \textbf{Human immunodeficiency virus (HIV)}
attacks the human immune
system by infecting the CD4+ cells (White Cells), which safeguards
the human body from diseases. HIV pushes hard the immune system to
its total failure within ten to
twelve years, leading to death if no treatment is sought.
If one is infected with HIV and have not taken any HIV drugs, then one
can expect three stages of illness as the infection progresses.
These three stages of HIV are:
\begin{enumerate}
    \item Acute HIV infection occurs shortly after HIV virus first enters the body. HIV
    infects the CD4+ cells and rapidly multiplies. Within a week or two,
    the infection takes hold and signs and symptoms of fever, headache,
    fatigue, and swollen lymph glands make many people feel like they
    are getting the flu. At this time, immune system is operating
    at full throttle, aggressively fighting the infection, and
    typically, after a few weeks, the flu-like symptoms will pass and
    person will start producing HIV antibodies~\cite{khan_academy}.

    \item Clinical latency, as its name suggests, is generally a time
    without symptoms. During this stage, your immune system and HIV
    have come
    into balance and the infection is partially controlled. However, HIV
    continues to slowly multiply and as time passes it steadily kills
    off the CD4+ cells making increasingly immunodeficient. A key
    point here is that during this time most people with HIV don't look
    or feel sick, and may not even know that they are infected.
    
    \item Acquired Immunodeficiency Virus (AIDS), it is the final stage of HIV in
    which less than $200$ CD4+ cells are
    left in the entire body. In this stage person survives around one to
    three years at most if and only if he/she gets treatment. AIDS
    starts when the immune system become so weak that it can no longer
    protect
    the person from infection by other organisms and diseases. One may
    experience any or all of the following symptoms: rapid weight loss,
    night sweats, extreme fatigue, swollen lymph glands, chronic
    diarrhea, sores in the mouth, anus or genitals, pneumonia, brown or
    purplish lesions on the skin or in the mouth, memory loss and
    depression (Figure~\ref{fig:figure_1}).
\end{enumerate}
\noindent
The HIV statistical estimates for $2002$ show that $4.69$ million South Africans are infected with HIV, and that by the year $2008$, half a million South Africans will die every year from AIDS-related causes~\cite{Shishana2004}. Recently, in $2019$ South Africa~\cite{HIV_africa} has made huge improvements in getting people to test for HIV in recent years and is now almost meeting the first of the $90-90-90$ targets, with $87$ percentage of people aware of their status. Where as, in latest surveillance report in Greece issued by the Hellenic Center For disease Control and Prevention (H.C.D.C.P) on October 
$2015$ reported that $15.109$ positive HIV infections. Among
these positive
HIV infections $3.782$ have already developed AIDS and
around $7.700$ are subject to antiretroviral therapy
(ART) and the deaths number due
to HIV infection has reached around $2.562$. In $2014$
report had showed huge number HIV positive people
acquired HIV due
sexual contact and injecting drug.
As the statistics suggest how important it is to find the cure for HIV infection, furthermore HIV infection has been studied in different articles~\cite{Callaway2002,Shishana2004,Verica2009,Burden2010,Charlotte2010,Kramer1999,Landi2008,Moysis2016,Rivadeneira2014,Rivadeneira2012,Nowak2001,Miguel2007,Mhawej2009,Mhawej2010,Xiaohua2006}.
\subsection{HIV dynamics}
\begin{figure}[ht]
    \centering
    \includegraphics[width=4.5in,height=2in]{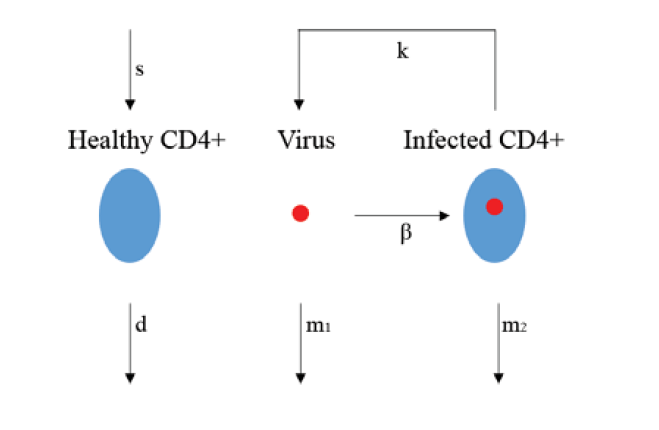}
    \caption{Interaction of Healthy CD4+ cells, Infected CD4+ Cells and Virus load, source:~\cite{Moysis2016}}
    \label{fig:figure_2}
\end{figure}

The simplest HIV dynamical model (\ref{appendix-section-3}) can be modelled with
the interaction of following cells and virus in the
human body: 
\begin{enumerate}
    \item Healthy CD4+ cells, $T(t)$
    \item Infected CD4+ cells, $T^\ast(t)$
    \item Viral Load (free virus or virions), $V(t)$
\end{enumerate}
\noindent
Healthy CD4+ cells are produced by the thymus at constant rate $s$
and also die at constant rate $d$. Healthy CD4+ cells are infected by virus
at the rate that is proportional to the  product of the number of
healthy cells and HIV Virus present in the body. The proportionality
constant known as the effectiveness of the infection is given by a
constant $\beta$. The infected cells come from healthy cells as they
get infected at the rate of $\beta  T(t)  V(t)$ and die at
constant rate $m_2$. The free virus particles, known as virions are
produced from infected CD4+ cells at a rate $k$ i.e. $kT^{*}(t)$ and
die at constant rate $m_1$ (Figure~\ref{fig:figure_2}).
\noindent
The interaction between the healthy CD4+ cells, infected CD4+ cells
and free virions (Virus Load) are modelled mathematically by the system of
nonlinear ordinary differential equations:
\begin{align} 
\frac{dT(t)}{dt}        &= s-dT(t)-\beta T(t)v(t), \label{HIV_System_11}\\
\frac{dT^\ast(t)}{dt}   &= \beta T(t)v(t)-m_2T^\ast(t), \label{HIV_System_12}\\
\frac{dV(t)}{dt}        &= kT^\ast(t)-m_1 v(t),\label{HIV_System_13}
\end{align}
\noindent
where $T(t)$ represents the number of healthy CD4+ cells,
$T^\ast(t)$ represents the number of infected CD4+ cells, and
$V(t)$ the number of virions, also known as the viral load. The
values for the different constant parameters used~\cite{Craig2005a,Miguel2007,Craig2004} in the system are given in the Table~\ref{table:ConstantsHIV}, 

\begin{table}[ht]
\centering 
\begin{tabular}{|p{6cm} c c|} 
\hline\hline 
Description & Symbol & Value \\ [0.5ex]
\hline \hline 
Time & $t$ & days \\ \hline
Death rate of uninfected T cells & $d$ & $0.02$ per day \\ 
\hline
Rate of virions produced & $k$ & $100$ counts per cell \\ \hline
Production rate of uninfected T cells & $s$ & $100$ mm$^{-3}$ per day \\ \hline
Ineffective  rate of virions & $\beta$ & $ 2.4\times10^{-5}$  mm$^{-3}$  per day \\ \hline
Death rate of virus & $ m_1 $ & $2.4$ per day \\ \hline
Death rate of infected T cells & $ m_2$ & $0.24$ per day\\\hline\hline
\end{tabular}
\caption{Values of the constant parameters used in system, (\ref{HIV_System_11})$-$(\ref{HIV_System_13}).} 
\label{table:ConstantsHIV} 
\end{table}
\subsection{Antiretroviral treatment}

The management of HIV/AIDS normally includes the use
of multiple antiretroviral drugs in an attempt to
control HIV infection. There are several classes of
antiretroviral agents that act on different stages of
the HIV life-cycle. The use of multiple drugs that act
on different viral targets is known as highly active
antiretroviral therapy (HAART). HAART decreases the
patient's total burden of HIV, maintains function of
the immune system, and prevents opportunistic
infections that often lead to death. HAART also
prevents the transmission of HIV between
serodiscordant same sex and opposite sex partners so
long as the HIV-positive partner maintains an
undetectable viral load.Treatment has been so
successful that in many parts of the world, HIV has
become a chronic condition in which progression to
AIDS is increasingly rare~\cite{WIKI_HIV}.

\noindent
HIV medicines work by preventing HIV from multiplying or by blocking infection of Healthy CD4+ cells known as reverse
transcriptase inhibitors (RTIs), while the
protease inhibitors (PIs) helps the immune system by preventing the production of
new virions, altogether RTIs and PIs help immune system to stabilise and suppress
the virus which results in decreased amount of virus in the bloodstream (viral
load). Although the
medicine does not get rid of HIV entirely, it gives CD4+ cells a
chance to recover so that they can fight off opportunistic
infections and cancers. If one don’t take ART then he/she is likely to die
within $12$ years from the time first got infected. On the other
hand, if one take ART then he/she can have a life expectancy equal to or
even higher than the general population.

\noindent
Taking into account the RTIs, $u_1(t)$, and PIs, $u_2(t)$, into the model  (\ref{HIV_System_11})$-$(\ref{HIV_System_13})
then the system takes the form as follows:
\begin{align} 
\frac{dT(t)}{dt}        &= s-dT(t)-(1-u_1(t))\beta T(t)v(t) \label{HIV_System_21}\\
\frac{dT^\ast(t)}{dt}   &= (1-u_1(t))\beta T(t)v(t)-m_2T^\ast(t) \label{HIV_System_22}\\
\frac{dV(t)}{dt}        &= (1-u_2(t))kT^\ast(t)-m_1 v(t)\label{HIV_System_23}
\end{align}

\noindent
The system $(\ref{HIV_System_21})-(\ref{HIV_System_23})$ is known as HIV dynamical model with treatment (Two input parameters), the terms
$(1 - u_1(t))$ and $(1 - u_2(t))$
represents the effectiveness of the RTIs and PIs,
respectively. If the values of $u_1(t) = 0$ and
$u_2(t) = 0$ then the treatment of HIV has not been
initiated and system $(\ref{HIV_System_21})-(\ref{HIV_System_23})$ becomes same as system~$(\ref{HIV_System_11})-(\ref{HIV_System_13})$.
If the values of $u_1(t) = 1$ and
$u_2(t) = 1$ then the treatment of HIV has worked
$100\%$ and system $(\ref{HIV_System_21})-(\ref{HIV_System_23})$ gets rid of free
virus and infected CD4+ cells resulting system in
Healthy CD4+ cells which is not attainable.

\noindent
The effectiveness of RTIs and PIs are further simplified by combining both RTIs
and PIs as single parameter $u(t)$ to administer the
drug. As seen in the clinical studies~\cite{Mhawej2009,Mhawej2010,Miguel2007}, hence it has been observed that combined treatment works better at preventing the production of
free virus (viral load) as compared to controlling spread of
infection. Keeping all these observation, the HIV dynamical model with treatment $(\ref{HIV_System_11})-(\ref{HIV_System_13})$ takes the form as follows:
\begin{align} 
\frac{dT(t)}{dt}        &= s-dT(t)-\beta T(t)v(t)\label{HIV_System_31} \\
\frac{dT^\ast(t)}{dt}   &= \beta T(t)v(t)-m_2T^\ast(t) \label{HIV_System_32}\\ 
\frac{dV(t)}{dt}        &= (1-u(t))kT^\ast(t)-m_1 v(t)\label{HIV_System_33}
\end{align}
\noindent
where the parameter, $u(t)$, represents the effectiveness
of the the combined treatment of RTIs and PIs.
\subsection{Conclusion}
The three HIV dynamical models have been proposed with
different clinical studies~\cite{Miguel2007,Mhawej2009,Mhawej2010}. Among these the first HIV dynamical model i.e. $(\ref{HIV_System_11})-(\ref{HIV_System_13})$ is considered, when there is no drug administered against HIV-Infection. HIV dynamical model with treatment including two input parameters RTIs and PIs is $(\ref{HIV_System_21})-(\ref{HIV_System_23})$ whereas with one (combined effect of RTIs and PIs) input parameter is $(\ref{HIV_System_31})-(\ref{HIV_System_33})$. Antiretroviral treatment is developed by highly active antiretroviral therapy (HAART), the aim of the antiretroviral treatment is to reduce the viral load to less than $50$ copies$/mm^3$. Moreover the HIV-Infection is incurable disease since last four decades, yet thanks to doctors and researchers that, who have managed to decrease the overall number of people developing to AIDS.

\section{Numerical simulation of HIV dynamics} \label{chapter2}
Numerical methods are used to approximate the solution of
equation when exact solution cannot be determined.
Numerical Methods help find the global behaviour of
ordinary differential equation. The HIV dynamical models discussed in Chapter~\ref{chapter1} are nonlinear and to approximate its solution the Runge-Kutta method of order $4$ is used~\cite{Burden2010}.

\subsection{Runge-Kutta method (RK4)}

\subsubsection{RK4 for $1^{\mathrm{st}}$ order ODE}
Consider a $1^{\mathrm{st}}$ order initial value problem of the form
\begin{align}\label{general_system_odes}
\frac{d u}{dt} &= f(t,u) \quad \mathrm{for} \quad a \leq t \leq b,\\
u(a)&=\alpha,
\end{align}
\noindent
The solution of system exist and is unique \cite{Burden2010}. However, in order to find the unknown functions $u(t)$ the classical Runge-Kutta $4^{\mathrm{th}}$ order method is used. The classical Runge-Kutta 4th order method used to solve the an ordinary differential equation (\ref{generalode}) is as follows,
\begin{align}
    &w_0 = \alpha \nonumber \\ 
    &k_1 = h f(t_i,w_i)  \nonumber \\
    &k_2 = h f(t_i + \frac{h}{2},w_i + \frac{k_1}{2}) \nonumber \\ 
    &k_3 = h f(t_i + \frac{h}{2},w_i + \frac{k_2}{2})\nonumber \\ 
    &k_4 = h f(t_i + h ,w_i + k_3)\nonumber \\ 
    &w_{i+1} = w_i + \frac{1}{6}(k_1 + 2k_2 + 2k_3 +k_4)
\end{align}
\noindent
Here $w_i$ represent the approximated solution of $u(t_i)$ on the specified
domain. 

\subsubsection{RK4 for system $1^{\mathrm{st}}$ order ODE}
Consider a $m$th order system of first order initial value problems of the form
\begin{align}\label{general_system_odes}
\frac{d u_{i}}{dt} &= f_{i}(t,u_{i}) \quad \mathrm{for} \quad a \leq t \leq b,\\
u_{i}(a)&=\alpha_{i},
\end{align}
\noindent
where $i= 1,2,\cdots, m$.
In order to solve the system of ODEs $(\ref{general_system_odes})-(\ref{general_system_odes})$, the Runge-Kutta method for system of $1^{\mathrm{th}}$ ordinary differential equations is briefly introduced here. The reader interested in detail is referred to \cite{Zill_ODEs_2012,Burden2010}.

\noindent
Let $N>0$ be an integer and set $h = (b-a)/N$ to partition the interval $[a,b]$ into $N$ sub-intervals with the mesh points $t_j = a +jh$, for each $j = 0,1,2,\cdots,N$. The notation $w_{ij}$, for each $j = 0,1,2,\cdots,N$ and $i = 1,2,\cdots,m$, is used to denote an approximated solution to $u_i(t_j)$. That is, $w_{ij}$ approximates the $i$th solution $u_i(t)$ of (\ref{general_system_odes}) at the $j$th mesh point $t_j$. For the initial conditions, set $w_{1,0} = \alpha_1, w_{2,0} = \alpha_2, \cdots,w_{m,0} = \alpha_m$ and assume that $w_{1,j}, w_{2,j}, \cdots, w_{m,j}$ are computed approximating $u_{1}(t_j), u_{2}(t_j), \cdots, u_{m}(t_j)$. In order to get $w_{1,j+1},w_{2,j+1},\cdots,w_{m,j+1}$, the following coefficient must be computed first, i.e.
\begin{align}
K_{1,i} &= hf_i(t_j,w_{1,j},w_{2,j},...,w_{m,j}), \\
K_{2,i} &= hf_i(t_j+0.5h,w_{1,j}+0.5k_{1,1},w_{2,j}+0.5k_{1,2},...,w_{m,j}+0.5k_{1,m}),\\
 K_{3,i} &= hf_i(t_j+0.5h,w_{1,j}+0.5k_{2,1},w_{2,j}+0.5k_{2,2},...,w_{m,j}+0.5k_{2,m}), \\
K_{4,i} &= hf_i(t_j+h,w_{1,j}+k_{3,1},w_{2,j}+k_{3,2},...,w_{m,j}+k_{3,m}),  
\end{align}
and
\begin{equation}
    w_{i,j+1} = w_{i,j} + \frac{1}{6}(k_{1,i} + 2k_{2,i} + 2k_{3,i} +k_{4,i}),
\end{equation}
Note that all the values $k_{1,1},k_{1,2},\cdots,k_{1,m}$ must be computed before any of the terms of the form $k_{2,i}$. In general, each $k_{p,1},k_{p,2},\cdots,k_{p,m}$ must be computed before any of the expressions $k_{p+1,i}$. 

\subsection{Numerical simulations of HIV dynamical models}
In this section, numerical simulation of HIV dynamical models proposed in the Chapter~\ref{chapter1} have been analyzed and discussed. The Runge-Kutta $4^\mathrm{th}$ order method is used to approximate the numerical solution of HIV dynamical
models \cite{Moysis2016,Craig2005a,Craig2004} under consideration. Algorithm of Runge-Kutta Method of order $4$ is implemented on MATLAB and the script are given in Appendix~\ref{appendix-b}).

%
\begin{figure}[htbp]
\begin{center}
\includegraphics[width=5in,height=1.6in]{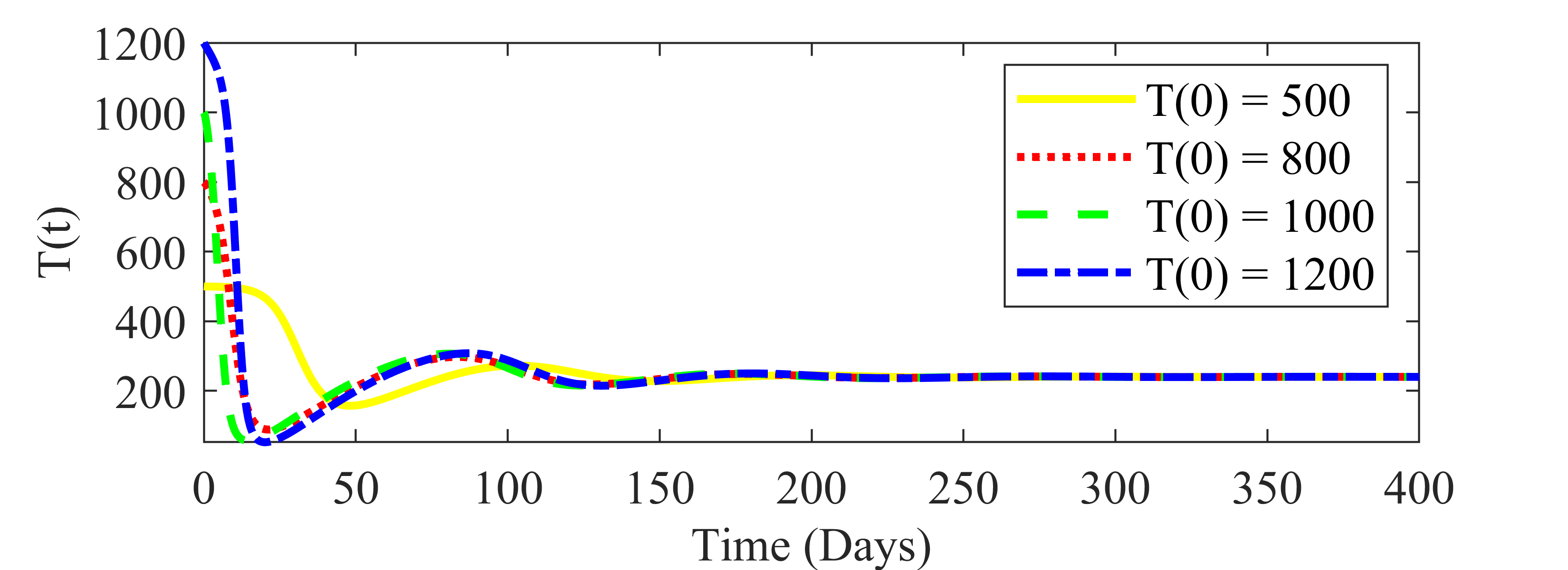} \\
\includegraphics[width=5in,height=1.6in]{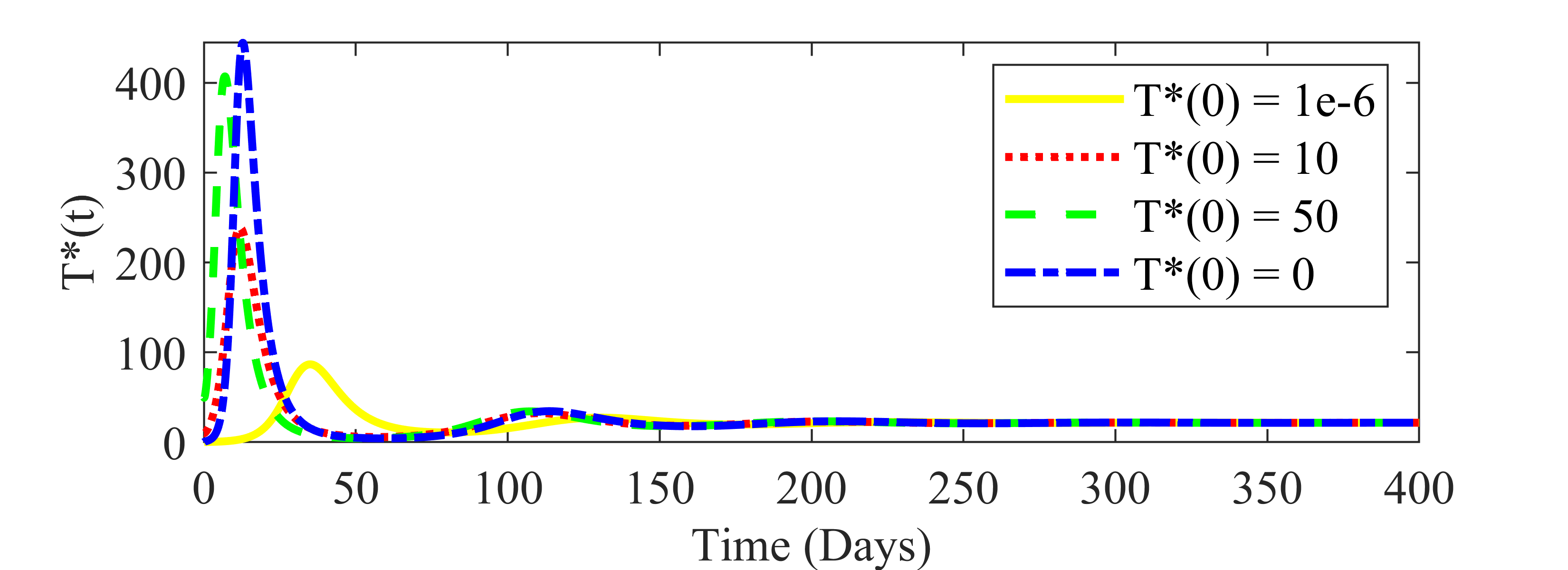}\\
\includegraphics[width=5in,height=1.6in]{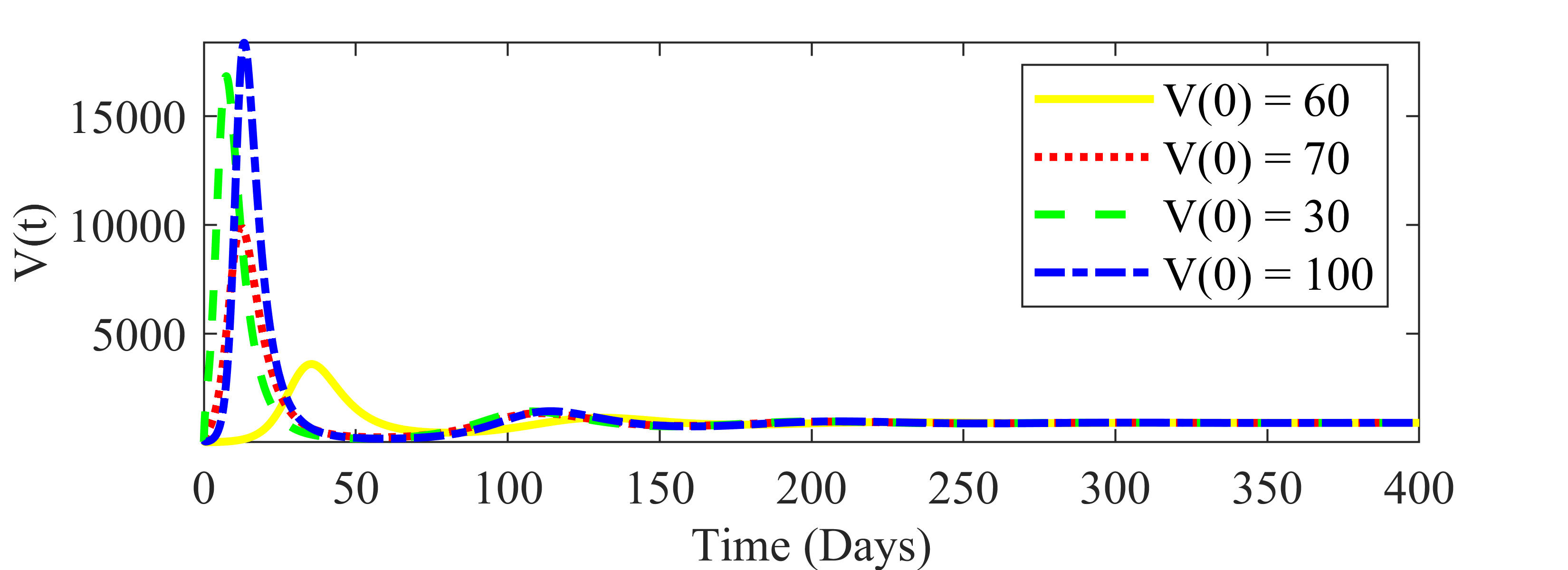} \\
\includegraphics[width=3.5in,height=2in]{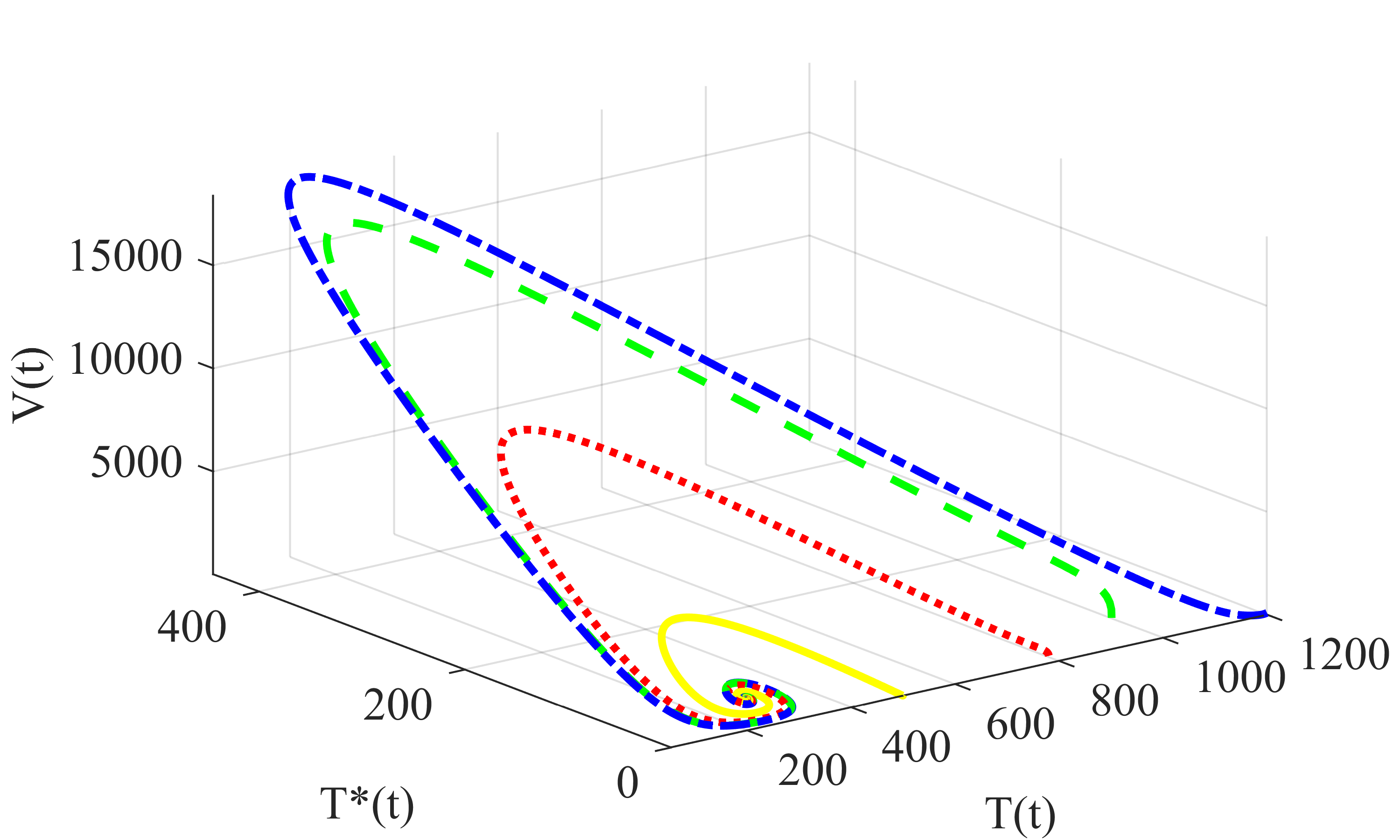} 
\end{center}
\caption{Simulations of HIV dynamics model without treatment, $(\ref{HIV_System_11})-(\ref{HIV_System_13})$.}
\label{pics:HIV_System_1_Figure}
\end{figure}
%

%
\begin{figure}[htbp]
\begin{center}
\includegraphics[width=5in,height=1.6in]{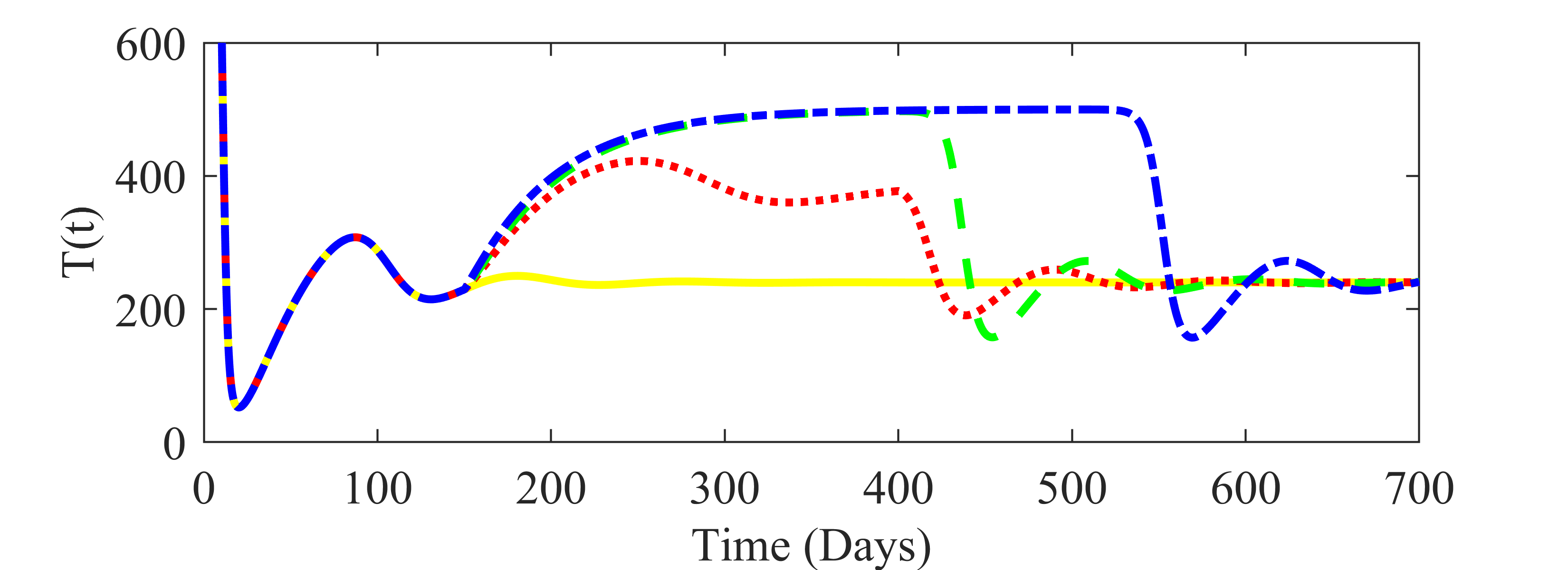} \\
\includegraphics[width=5in,height=1.6in]{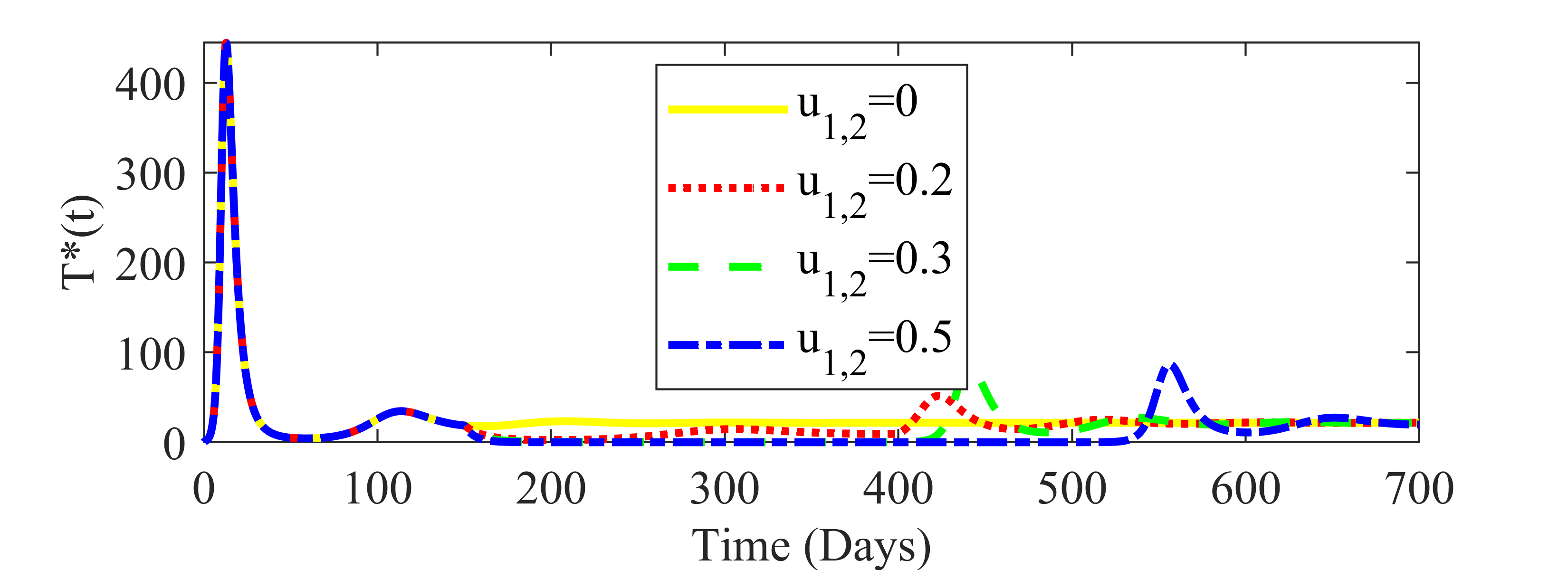}\\
\includegraphics[width=5in,height=1.6in]{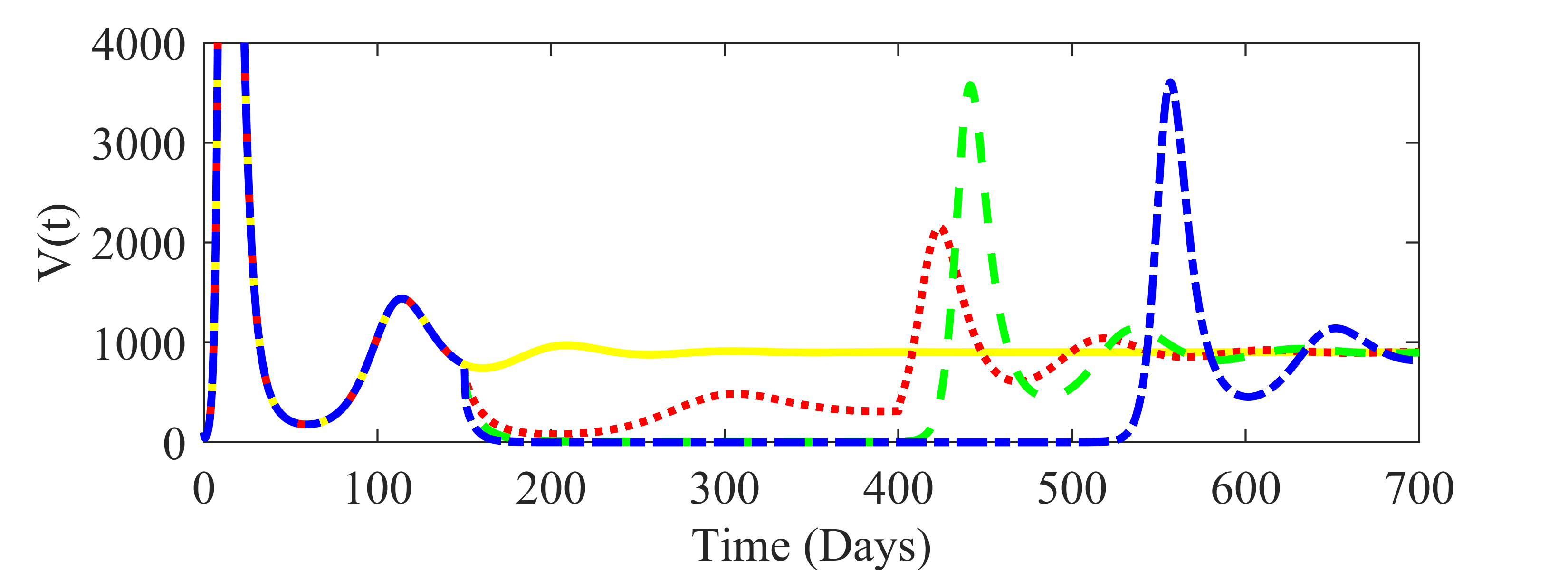} \\
\includegraphics[width=3.5in,height=2in]{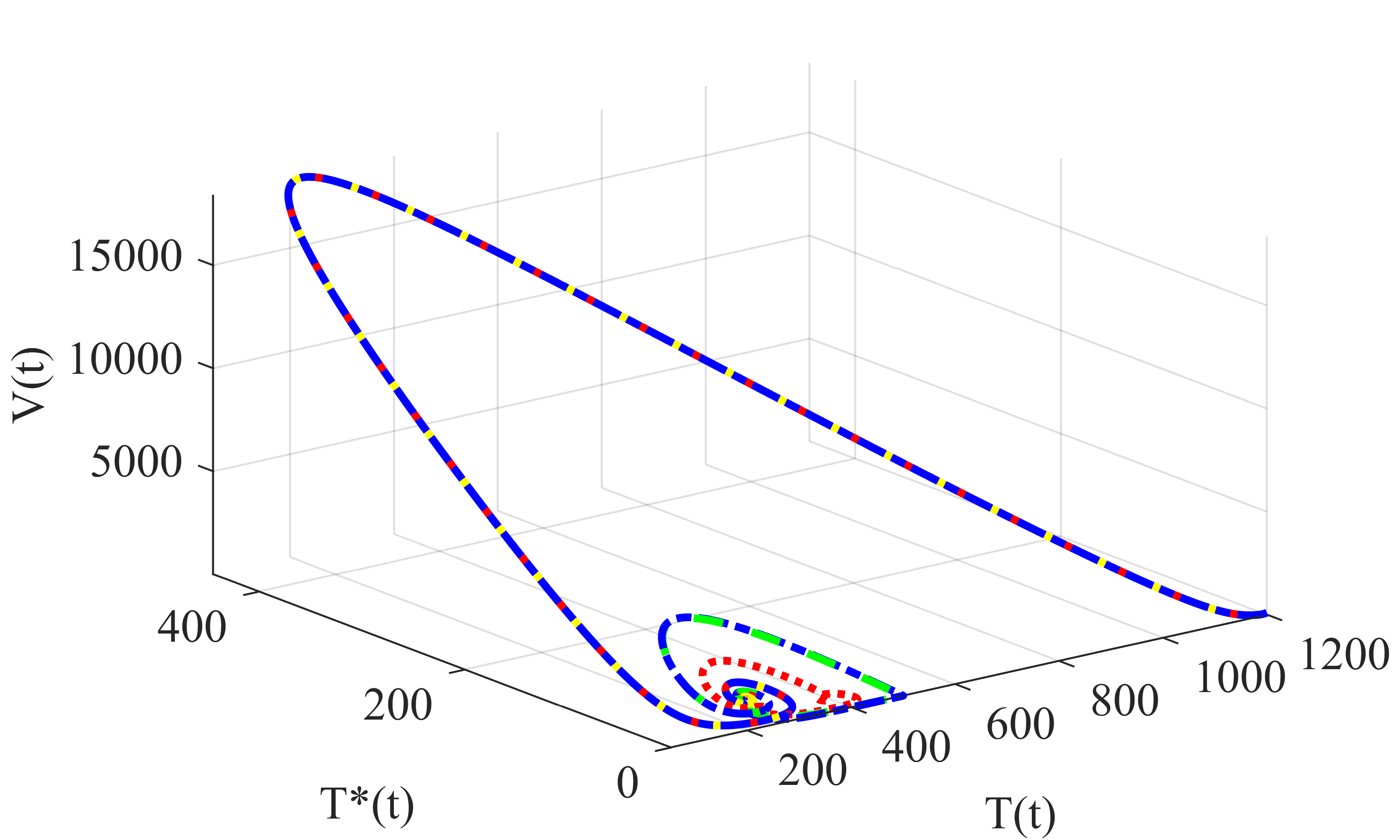} 
\end{center}
\caption{Simulations of HIV dynamics model with treatment, $(\ref{HIV_System_21}) - (\ref{HIV_System_23})$.}
\label{pics:HIV_System_2_Figure}
\end{figure}
%
%
\begin{figure}[htbp]
    \begin{center}
        \includegraphics[width=5in,height=1.6in]{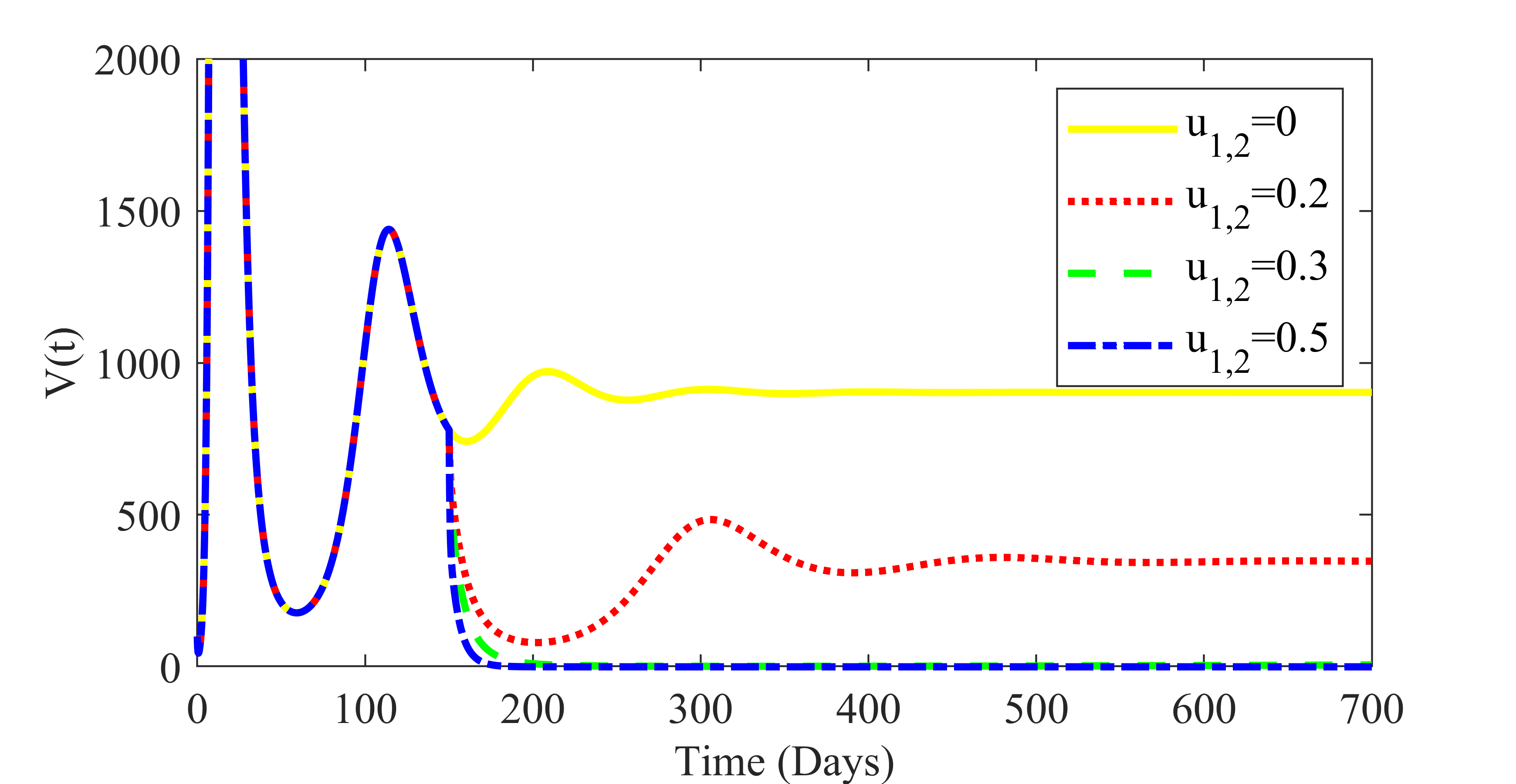}
    \end{center}
\caption{Simulations of Viral Load in HIV Dynamics model with treatment, $(\ref{HIV_System_21}) - (\ref{HIV_System_23})$. Treatment initiated at $150^{\mathrm{th}}$ day and is taken without any break or interruption.}
\label{pics:HIV_System_2_Figure_Continuous_Treatment}
\end{figure}
%
\subsubsection{HIV dynamical model without treatment}

In this section, the HIV dynamical model $(\ref{HIV_System_11}) - (\ref{HIV_System_13})$ in which progression of CD4+ healthy cells, infected CD4+ cells and viral load are analyzed, numerically, as the time passes (Figure~\ref{pics:HIV_System_1_Figure}). 

\noindent
Initially, when someone gets infected by HIV virus then in the first three weeks there is sudden decrease in the healthy CD4+ cells, and healthy CD4+ cells level decreases to $200$ cells (Figure~\ref{pics:HIV_System_1_Figure}$a$). The infected CD4+ cells level increases to $500$ cells (Figure~\ref{pics:HIV_System_1_Figure}$b$). Whereas the HIV virus level reaches above 15000 viral load (Figure~\ref{pics:HIV_System_1_Figure}$c$). This is called acute HIV-Infection stage, and for few people it may last up to $6$ months. 

\noindent
After $7$ to $12$ weeks, the immune system starts suppressing the HIV virus due to its abundance presence in the blood. And after $150$ days, there is a balance between healthy CD4+ cells and viral load, although the HIV virus is still active but its reproductive cycle is repressed by the immune system. This stage is called clinical latency stage, it may last up to $10$ years.

\noindent
A typical progression for the disease is shown in (Figure~\ref{pics:HIV_System_1_Figure}$d$). It is clear that after initial
infection as in the $2^{\mathrm{nd}}$ paragraph of this section, there is a rise in the infected CD4+ cells and after the reaction of the immune system,
the system is stabilized and the patient enters the clinical latency stage. It is observed that, no matter what the initial condition of patient's (Healthy CD4+ cells, Infected CD4+ cells and Virus Load), trajectories of HIV dynamical model $(\ref{HIV_System_11})-(\ref{HIV_System_13})$ always reaches the same equilibrium point i.e. $(T=240,\; T^\ast=21.6667,\; V=902.778)$. In this model spread of virus and has been observed and no treatment is assumed. Initial conditions used to simulate the model has been mentioned in the graph.

\begin{figure}[hp]
\begin{center}
\includegraphics[width=5in,height=1.6in]{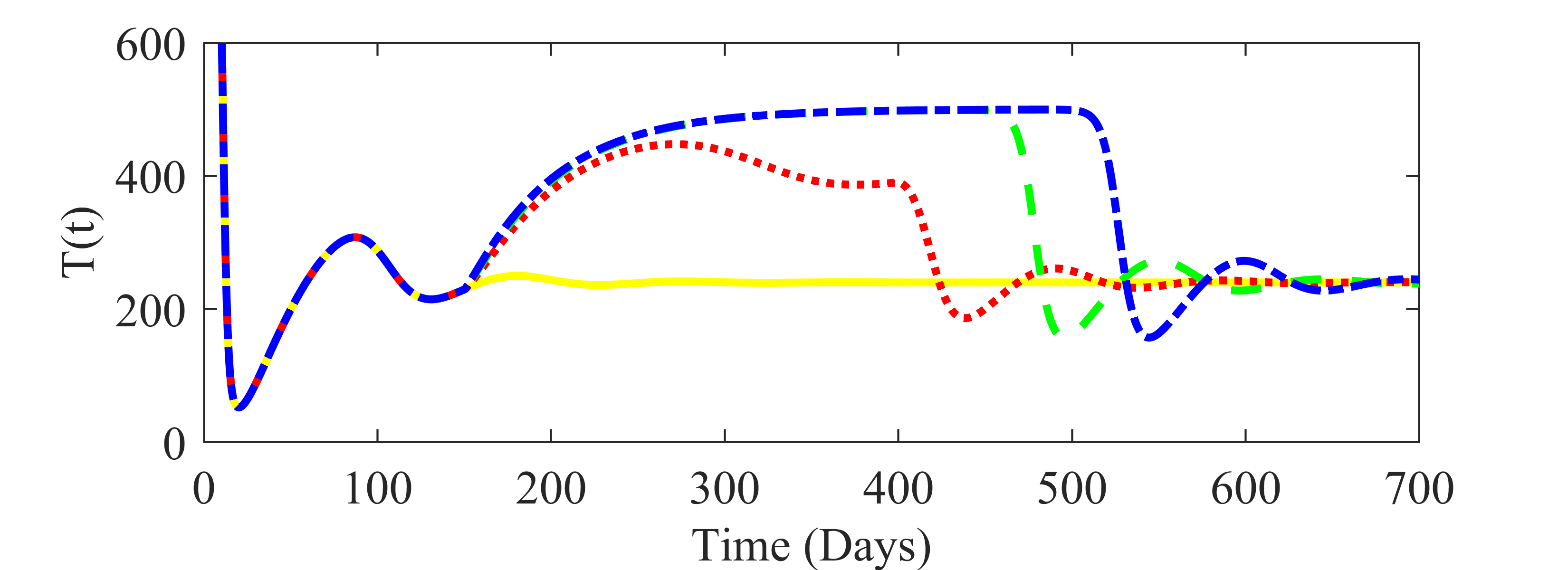} \\
\includegraphics[width=5in,height=1.6in]{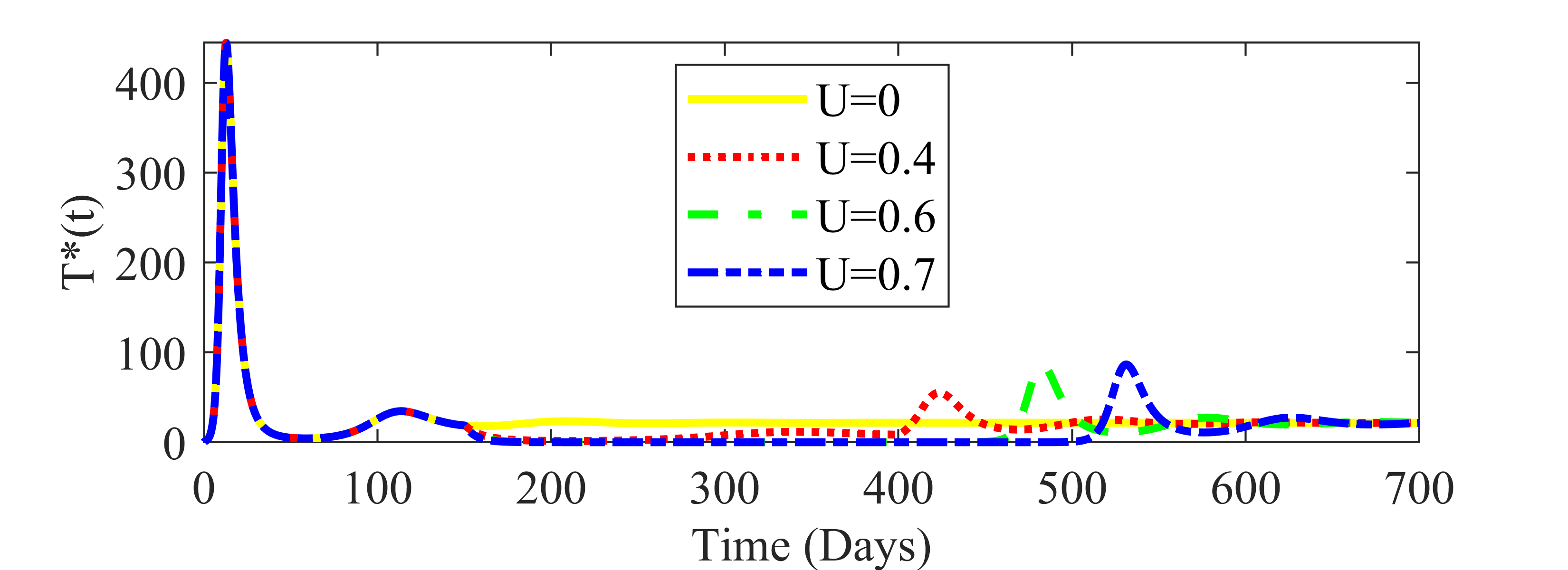}\\
\includegraphics[width=5in,height=1.6in]{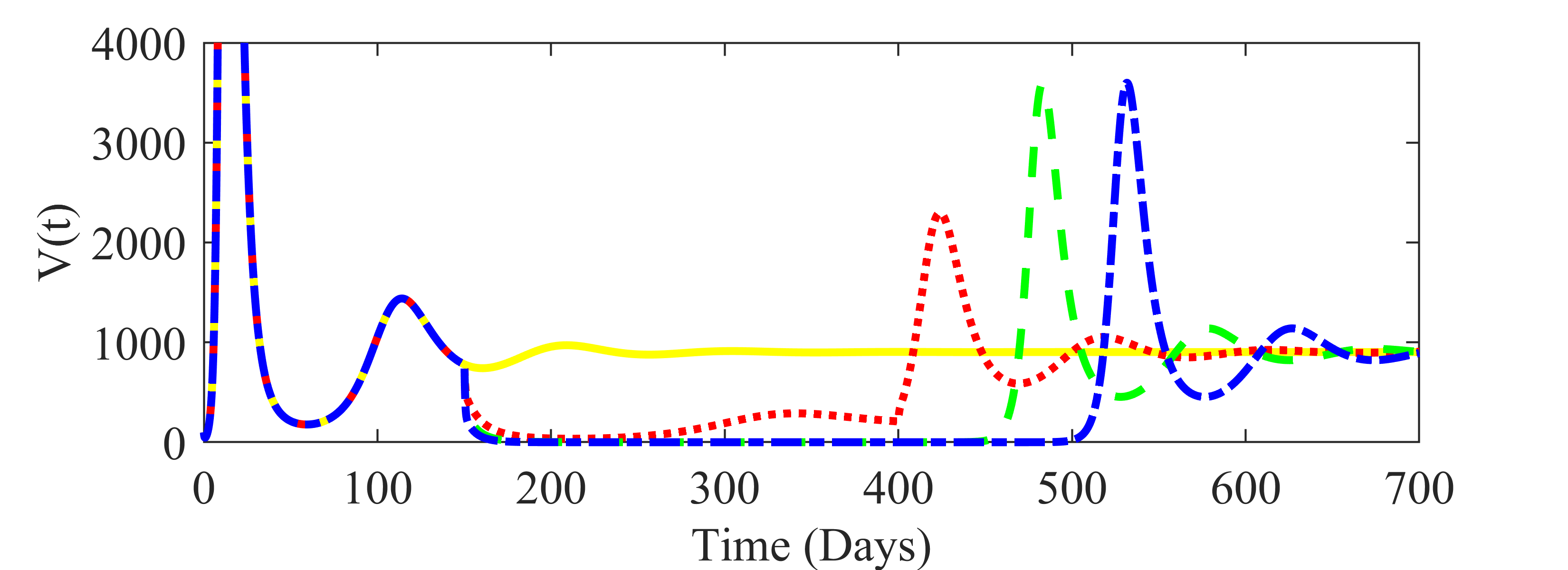} \\
\includegraphics[width=3.5in,height=2in]{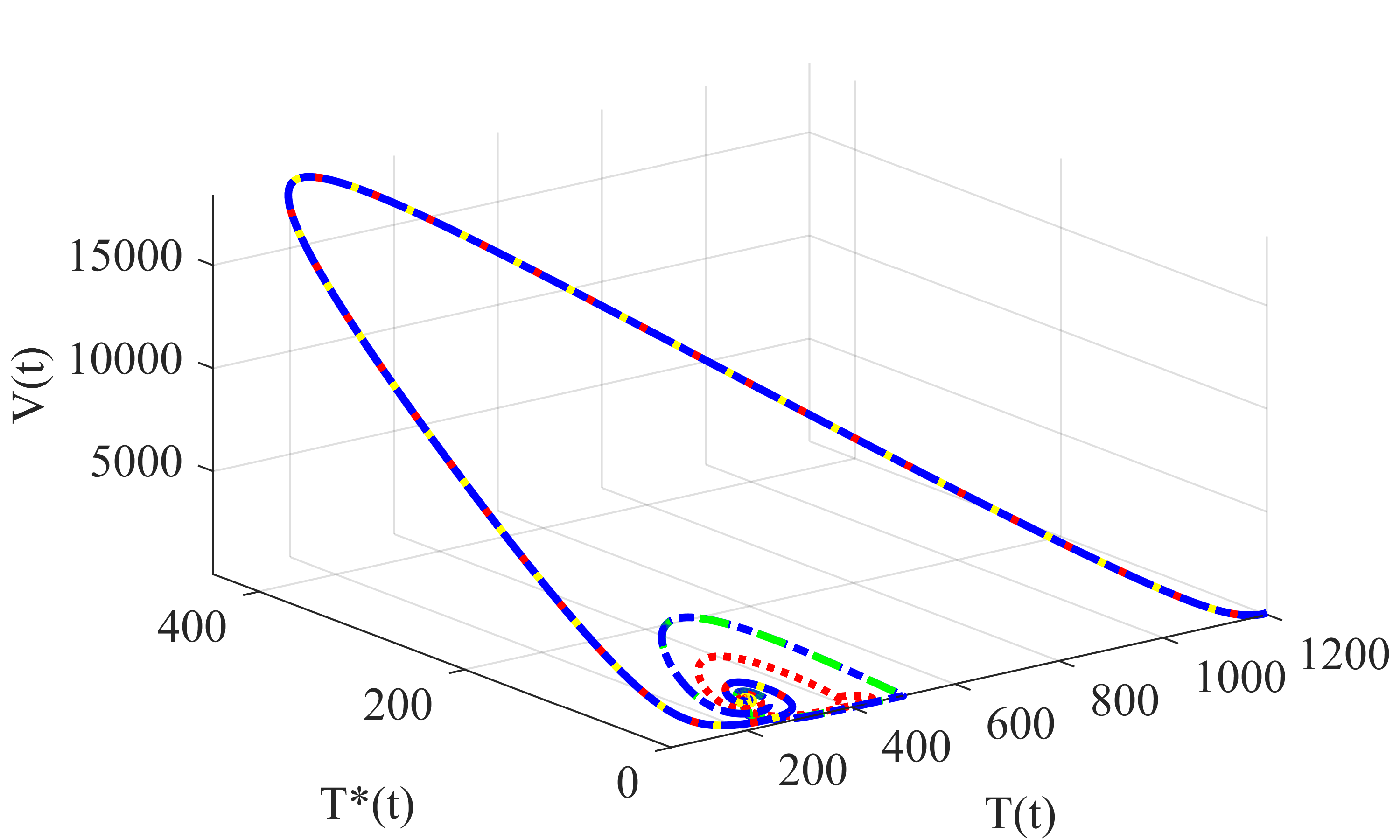} 
\end{center}
\caption{Simulations of HIV Dynamics model with treatment,  $(\ref{HIV_System_31}) - (\ref{HIV_System_33})$.}
\label{pics:HIV_System_3_Figure}
\end{figure}

\subsubsection{HIV dynamical model With treatment}

\paragraph{With two input parameter}

In this section, HIV dynamical model $(\ref{HIV_System_21})-(\ref{HIV_System_23})$ in which the additional treatment parameters $i)$ Reverse Transciptase Inhibitors, $u_1(t)$ and $ii)$ Protease Inhibitors, $u_2(t)$, are incorporated.

\noindent
The numerical simulations are performed by keeping the fixed initial conditions $T(0)=1200$, $T^\ast(0)=0$ and $V(0)=100$ of Healthy CD4+ cells, Infected CD4+ cells and Virus load, respectively. The treatment parameters $u_1$ and 
$u_2$ represents, how much, does the treatment succeed in combating the HIV-Infection.
If $u_1 = 0$ and $u_2 = 0$ then no therapy is considered. The solid line in the Figure \ref{pics:HIV_System_2_Figure} represents the no therapy condition and the model works same as that of HIV model without treatment $(\ref{HIV_System_11})-(\ref{HIV_System_13})$. Here, the low and medium low drug dosage is represented by setting $u_{1,2} = 0.2$ and  $u_{1,2} = 0.3$, respectively. Moreover, $u_{1,2} = 0.5$ is considered for medium drug dosage.

\noindent
It is seen in Figure \ref{pics:HIV_System_2_Figure} that, at $400^{\mathrm{th}}$ day treatment has been terminated which results in virus load rise again. Even though the treatment is terminated, it is seen that virus load stay low for certain number of days depending on the effectiveness of the treatment. 
In order to see the effectiveness of treatment, the model has been initiated at $150^{\mathrm{th}}$ day, and before $150^{\mathrm{th}}$ day it is assumed that the person has not gone under any treatment. Virus load and infected CD4+ cells get suppressed after $150^{\mathrm{th}}$ day, as witnessed, in the simulation (Figure~\ref{pics:HIV_System_2_Figure_Continuous_Treatment}).

\noindent
 In order to suppress the virus load as well as the infected CD4+ Cells it is necessary that treatment should be taking continuous, once it is initiated. This is witnessed in the Figure \ref{pics:HIV_System_2_Figure_Continuous_Treatment}. 

\paragraph{With one input parameter}
In this model the treatment effect of reverse transcriptase inhibitors and protease inhibitors has been combined as a single parameter $U(t)$.
It is already assessed clinically in \cite{Miguel2007,Mhawej2009,Mhawej2010} that effect of RTIs and PIs should be combined to combat the HIV Virus.

\noindent
In the simulation, $U=0$ represent viral load under no therapy, $U=0.4$ represent viral load for low drug dosage, $U=0.6$ represent viral load for medium drug dosage and $U=0.7$ viral load for high drug dosage.
It is observed (Figure \ref{pics:HIV_System_3_Figure}) that the viral load bounce back immediately making small oscillations. Initially the person is not treated and model works same as that of $(\ref{HIV_System_11})-(\ref{HIV_System_13})$ and at $150^{\mathrm{th}}$ day treatment has been initiated with the model $(\ref{HIV_System_31}) - (\ref{HIV_System_33})$ with different drug dosages $U =0, U=0.4, U=0.6$ and  $U=0.7$.
The treatment is then terminated at $400^{\mathrm{th}}$ day which results in rise of viral load and infected CD4+ cells after certain interval of time, that depends on how much does the drug is effective. Overall, the result of this model suggests that if combined drug is added (as witnessed in this model), then the results for both the models $(\ref{HIV_System_21})-(\ref{HIV_System_23})$ and $(\ref{HIV_System_31})-(\ref{HIV_System_33})$ are same. However, the benefit of this model is that it depends only on one parameter, $U(t)$. And one of the drawback is that viral load comes back faster as compare to HIV dynamical model with treatment (two input parameters)

\subsection{Conclusion}
It is concluded that when HIV dynamical model with no treatment is simulated, the trajectories converge to equilibrium point which suggests that the person has entered in the clinical latency stage. The trajectories of HIV dynamical model with treatment (two input parameter) suggests that viral load is suppressed very efficiently if the treatment is taken without any interruption, and if the treatment is terminated, still the virus remain low for certain days. whereas the HIV dynamical with treatment (one input parameter) shows results better although in this model only production of new virions is minimised up to $70\%$ and viral load is reduced to 50 copies per millimeter cube, but if the treatment is terminated the viral load bounce back faster as compared to two input parameter model. Finally, the benefit of using RK4 method of order $4$ is that the global solution of nonlinear system of ODEs is found subject to given initial conditions.

\section{Linearization}\label{chapter3}
This chapter present the linearization of the dynamical system, which is used to analyze behaviour around the equilibrium or critical points of the given nonlinear system. 

\subsection{Linearization Technique}
A general first order autonomous $3 \times 3$ system of ODEs has the following form:
\begin{align}
\frac{dx(t)}{dt} &= f_1(x(t),y(t),z(t)) \label{gsystem_11} \\ 
\frac{dy(t)}{dt} &= f_2(x(t),y(t),z(t)) \label{gsystem_12}\\
\frac{dz(t)}{dt} &= f_3(x(t),y(t),z(t))\label{gsystem_13}
\end{align}
\noindent
where $x(t), \; y(t)$ and $z(t)$ are the solution of system $(\ref{gsystem_11})-(\ref{gsystem_13})$.
The critical point, say $( x_0,y_0,z_0)$, of the system $(\ref{gsystem_11})-(\ref{gsystem_13})$ are those,  where ${dx(t)}/{dt}=0$ and ${dy(t)}/{dt}=0$ and ${dz(t)}/{dt}=0$, which satisfies the following
\begin{align} \label{gsystem_2}
f_1(x_0,y_0,z_0)&=0 \\ 
f_2(x_0,y_0,z_0)&=0 \\
f_3(x_0,y_0,z_0)&=0 
\end{align}
\noindent
Equivalently, these are the equilibrium solution. 
Now with the help of change of variables i.e.
\begin{align}
     u = x - x_0 &=\Delta x \\
     v = y - y_0 &=\Delta y\\
     w = z - z_0 &=\Delta z 
\end{align}
\noindent
Above change of variables will put origin of uvw-plane at $(x_0,y_0,z_0)$ and similarly
\begin{align}
     \frac{du(t)}{dt} &=\frac{dx(t)}{dt}\\
     \frac{dv(t)}{dt} &=\frac{dy(t)}{dt}\\
     \frac{dw(t)}{dt} &=\frac{dz(t)}{dt}
\end{align}
\noindent
Using \(x = x_0 + u\), \(y = y_0 + v\) and \(z = z_0 + w\) in $(\ref{})$ yields
\begin{align} \label{gsystem_3}
f_1(u+x_0,v+y_0,w+z_0) &= f_1(x_0,y_0,z_0) +  \frac{\partial}{\partial
x}f_1(x_0,y_0,z_0)(\Delta x) + \\ \nonumber & \frac{\partial}{\partial
y}f_1(x_0,y_0,z_0)(\Delta y) +\frac{\partial}{\partial
z}f_1(x_0,y_0,z_0)(\Delta z) \\
f_2(u+x_0,v+y_0,w+z_0) &= f_2(x_0,y_0,z_0) + \frac{\partial}{\partial
x}f_2(x_0,y_0,z_0)(\Delta x) + \\ \nonumber & \frac{\partial}{\partial
y}f_2(x_0,y_0,z_0)(\Delta y) +\frac{\partial}{\partial
z}f_2(x_0,y_0,z_0)(\Delta z) \\
f_3(u+x_0,v+y_0,w+z_0)&= f_3(x_0,y_0,z_0) + \frac{\partial}{\partial
x}f_3(x_0,y_0,z_0)(\Delta x) + \\ \nonumber & \frac{\partial}{\partial
y}f_3(x_0,y_0,z_0)(\Delta y) +\frac{\partial}{\partial
z}f_3(x_0,y_0,z_0)(\Delta z) 
\end{align}

\noindent The above set of equations can be rewritten in the matrix form
\begin{align}
&\begin{bmatrix}
f_1(u+x_0,v+y_0,w+z_0) \\
f_2(u+x_0,v+y_0,w+z_0)\\
f_3(u+x_0,v+y_0,w+z_0)
\end{bmatrix}
 = 
\begin{bmatrix}
f_1(x_0,y_0,z_0) \\
f_2(x_0,y_0,z_0)\\
f_3(x_0,y_0,z_0)
\end{bmatrix}
+ \\ \nonumber &
\begin{bmatrix}
\frac{\partial}{\partial x}f_1(x_0,y_0,z_0) &\frac{\partial}{\partial y}f_1(x_0,y_0,z_0) &\frac{\partial}{\partial z}f_1(x_0,y_0,z_0) \\
\frac{\partial}{\partial x}f_2(x_0,y_0,z_0) &\frac{\partial}{\partial y}f_2(x_0,y_0,z_0) &\frac{\partial}{\partial z}f_2(x_0,y_0,z_0) \\
\frac{\partial}{\partial x}f_3(x_0,y_0,z_0) &\frac{\partial}{\partial y}f_3(x_0,y_0,z_0) &\frac{\partial}{\partial z}f_3(x_0,y_0,z_0) \\
\end{bmatrix}
\begin{bmatrix}
\Delta x\\
\Delta y\\
\Delta z
\end{bmatrix}
\end{align}
\noindent
Now as \(f_i(x_0,y_0,z_0) = 0 \) for \(i=1,2,3\) and \(\Delta x = u\) , \(\Delta y = v\) and \(\Delta z = w\),
the system in matrix form follows
$$
\begin{bmatrix}
\frac{du(t)}{dt}\\
\frac{dv(t)}{dt}\\
\frac{dw(t)}{dt}
\end{bmatrix}
=
\begin{bmatrix}
\frac{\partial}{\partial x}f_1(x_0,y_0,z_0) &\frac{\partial}{\partial y}f_1(x_0,y_0,z_0) &\frac{\partial}{\partial z}f_1(x_0,y_0,z_0) \\
\frac{\partial}{\partial x}f_2(x_0,y_0,z_0) &\frac{\partial}{\partial y}f_2(x_0,y_0,z_0) &\frac{\partial}{\partial z}f_2(x_0,y_0,z_0) \\
\frac{\partial}{\partial x}f_3(x_0,y_0,z_0) &\frac{\partial}{\partial y}f_3(x_0,y_0,z_0) &\frac{\partial}{\partial z}f_3(x_0,y_0,z_0) \\
\end{bmatrix}
\begin{bmatrix}
u(t)\\
v(t)\\
w(t)
\end{bmatrix}
$$
\noindent
Here, the Jacobian matrix $J$ is as follows
\begin{align}
    J(x_0,y_0,z_0) =
\begin{bmatrix}
\frac{\partial}{\partial x}f_1(x_0,y_0,z_0) &\frac{\partial}{\partial y}f_1(x_0,y_0,z_0) &\frac{\partial}{\partial z}f_1(x_0,y_0,z_0) \\
\frac{\partial}{\partial x}f_2(x_0,y_0,z_0) &\frac{\partial}{\partial y}f_2(x_0,y_0,z_0) &\frac{\partial}{\partial z}f_2(x_0,y_0,z_0) \\
\frac{\partial}{\partial x}f_3(x_0,y_0,z_0) &\frac{\partial}{\partial y}f_3(x_0,y_0,z_0) &\frac{\partial}{\partial z}f_3(x_0,y_0,z_0) \\
\end{bmatrix}
\end{align}
\noindent
and the corresponding linearized system is given by
\begin{equation}\label{linearizedsystem}
    \begin{bmatrix}
    \frac{du(t)}{dt}\\
    \frac{dv(t)}{dt}\\
    \frac{dw(t)}{dt}
    \end{bmatrix}
    =
    J(x_0,y_0,z_0)
    \begin{bmatrix}
    u(t)\\
    v(t)\\
    w(t)
    \end{bmatrix}
\end{equation}
\noindent
This is a linearized system of the given nonlinear ordinary differential equations with coefficient matrix $J(x_0,y_0,z_0)$. It is called the linearization of the system around the critical point $(x_0,y_0,z_0)$. Hence, the solution of the linearized system (\ref{linearizedsystem}) is written as.
\begin{equation}
    X(t) = c_1 V_1 e^{\lambda_1 t} + c_2 V_2 e^{\lambda _2 t} + c_3 V_3 e^{\lambda _3 t}.
\end{equation}
\noindent where $c_1, \; c_2$ and $c_3$ are any real constant and $\lambda_1, \; \lambda_2$ and $\lambda_3$ are the eigenvalues of Jacobian matrix (i.e. $J(x_0,y_0,z_0)$) and $V_1, \; V_2$ and $V_3$ are vector associated eigenvectors.

\subsection{Linearization of HIV dynamical model without treatment}
Let us now linearize the nonlinear system of the HIV model 
\begin{align} \label{system_1}
\frac{dT(t)}{dt} &= f_1(T(t),T^\ast(t),V(t)) \\
\frac{dT^\ast(t)}{dt} &= f_2(T(t),T^\ast(t),V(t))\\
\frac{dv(t)}{dt} &= f_3(T(t),T^\ast(t),V(t))
\end{align}
where $f_1, \; f_2$ and $f_3$ are as follows:
\begin{align} \label{system_2}
f_1(T(t),T^\ast(t),V(t)) &= s-dT(t)-\beta T(t)v(t) \\
f_2(T(t),T^\ast(t),V(t))&= \beta T(t)v(t)-m_2T^\ast(t) \\
f_3(T(t),T^\ast(t),V(t))&= kT^\ast(t)-m_1 v(t)
\end{align}

\noindent
In order to linearize the HIV dynamical model $(\ref{HIV_System_11})-(\ref{HIV_System_13})$, one need to first find the critical points of the given system. The first critical point of the system is 
\begin{equation} \label{critical_1}
\left( \frac{s}{d} , 0 , 0 \right)
\end{equation}
\noindent
This critical point (\ref{critical_1}) indicates that neither a person have any viral load nor it has infected CD4+ cells. Hence, it is concluded that a person is healthy and may be happy as well. The second critical point of the system is
\begin{equation} \label{critical_2}
\left( \frac{m_1 m_2}{k\beta} , \frac{s}{m_2} - \frac{d m_1}{\beta k} , \frac{k s}{m_1 m_2} - \frac{d}{\beta} \right).
\end{equation}
%

\begin{figure}[htbp]
\begin{center}
\includegraphics[width=5.5in,height=1.6in]{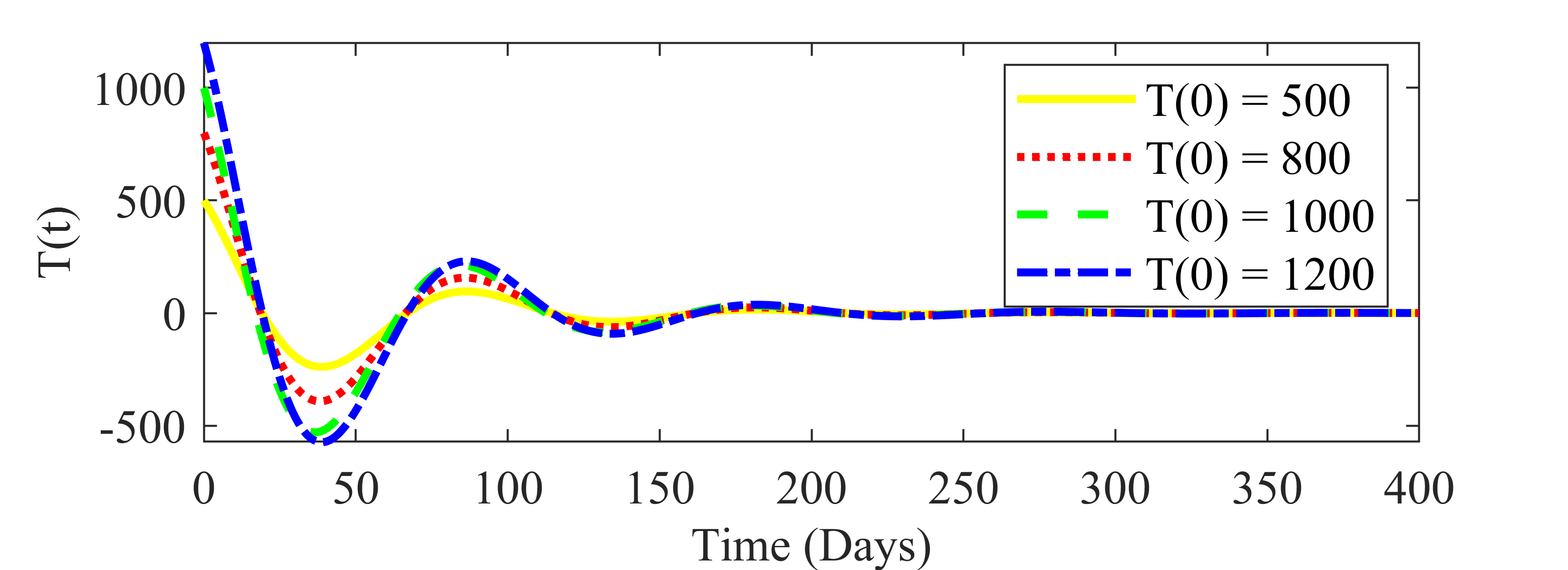} \\
\includegraphics[width=5.5in,height=1.6in]{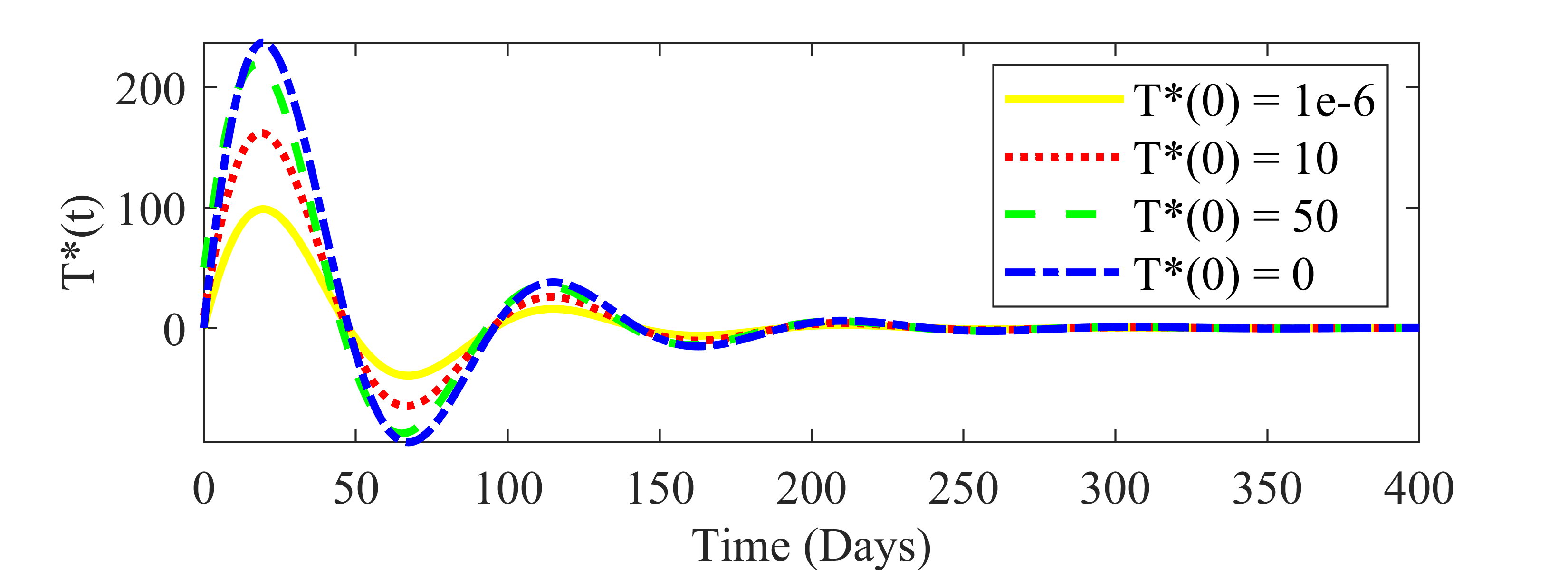}\\
\includegraphics[width=5.5in,height=1.6in]{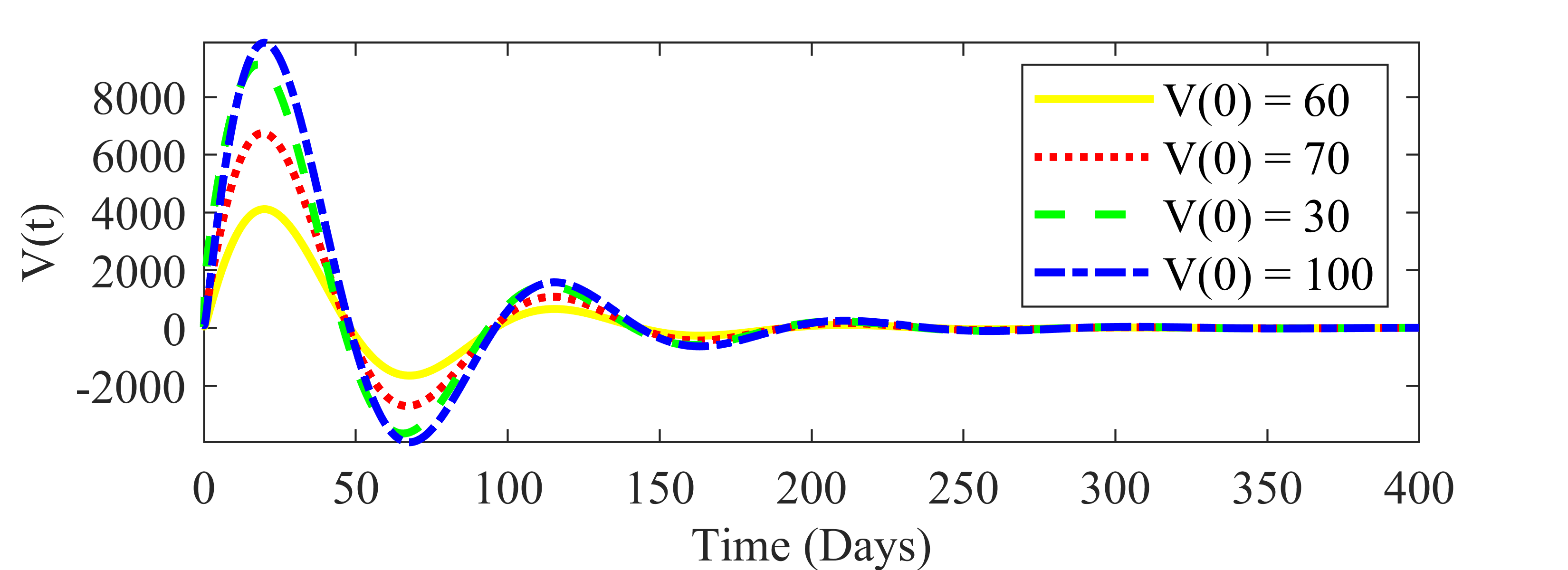} \\
\includegraphics[width=3.5in,height=2in]{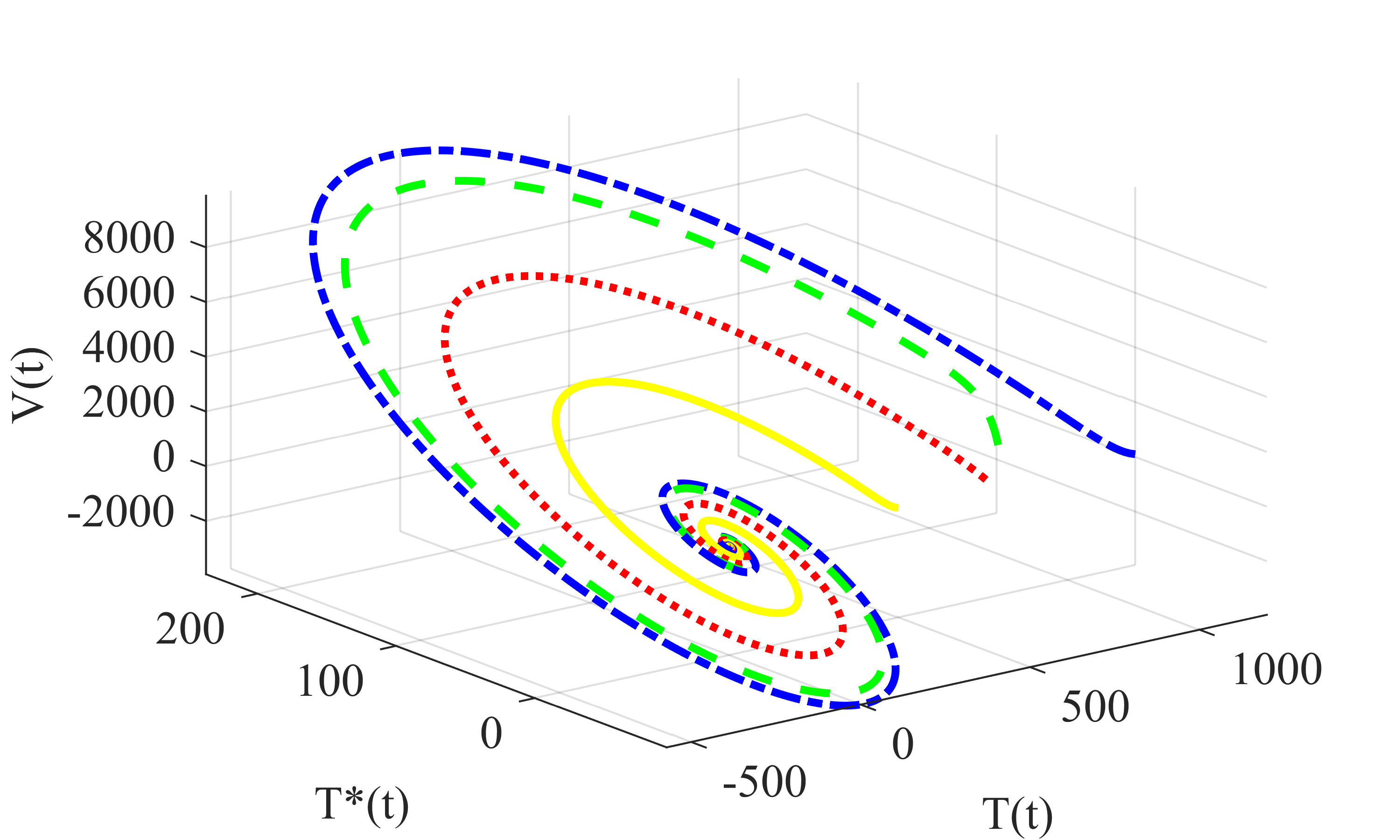} 
\end{center}
\caption{Phase Plane and Portrait of linearized HIV dynamics model without treatment, $(\ref{linearized_model_one})$.}
\label{pics:HIV_System_1_Figure_linearized}
\end{figure}

%
\noindent The critical point~(\ref{critical_2}) shows that a person is passing through the clinical stage of the HIV~(\cite{Craig2005a}). Now to see how does the system of ODEs of HIV dynamical model~$(\ref{HIV_System_11})-(\ref{HIV_System_13})$ behave around critical point~\ref{critical_2}. By using the values given in Table~\ref{table:ConstantsHIV} the critical point exact numerical value is computed. i.e.
\begin{equation}
\left(\frac{m_1 m_2}{k\beta} , \frac{s}{m_2} - \frac{d m_1}{\beta k} , \frac{k s}{m_1 m_2} - \frac{d}{\beta}\right) =   \left(240 , 21.6667 , 902.778 \right)  
\end{equation}
\noindent
Now by using the theory of linearization \cite{Zill_ODEs_2012,Nagle2011}, Jacobian of HIV dynamical model $(\ref{HIV_System_11})-(\ref{HIV_System_13})$ is given by
\begin{align}
J \left(240 , 21.6667 , 902.778 \right)  = 
\begin{bmatrix}
-0.0417 & 0   & -0.0058 \\
0.0217   & -0.24 & 0.0058 \\
0          & 100   &  -2.4 \\
\end{bmatrix}
\end{align}
\noindent
This implies that the linearized problem of the nonlinear HIV system under discussion is
\begin{align}\label{linearized_model_one}
\begin{bmatrix}
u^{'}(t)\\
v^{'}(t)\\
w^{'}(t)
\end{bmatrix} =\begin{bmatrix}
-0.0417 & 0   & -0.0058 \\
0.0217   & -0.24 & 0.0058 \\
0          & 100   &  -2.4 \\
\end{bmatrix}
\begin{bmatrix}
u(t)\\
v(t)\\
w(t)
\end{bmatrix}
\end{align}
\noindent
The eigenvalues of Jacobian \cite{Craig2004,Craig2005a} that are
\begin{align} 
    \lambda _1 & = -2.64334 \label{eigenvalue1} \\
    \lambda _2  & = -0.0191783 + 0.0658064\iota \label{eigenvalue2} \\
    \lambda _3 &= -0.0191783 - 0.0658064\iota  \label{eigenvalue3}
\end{align}
\noindent
The real part of all the eigenvalues $(\ref{eigenvalue1})-(\ref{eigenvalue3})$ are negative, this implies that corresponding critical point (\ref{critical_2}) on which the system has been linearized is an asymptotically stable node and also real part of eigenvalues are nonzero which implies that the equilibrium point is hyperbolic point. Furthermore, with the Hartman–Grobman theorem implies that when a critical point is hyperbolic point then the linearized system around the hyperbolic point preserves the qualitative behaviour of the corresponding nonlinear system i.e. HIV dynamical model $(\ref{HIV_System_11})-(\ref{HIV_System_13})$. The eigenvectors corresponding to eigenvalues are

\begin{equation}
        V_1 = 
    \begin{bmatrix}
        0.0022 \\
        -0.0024 \\
        1  
    \end{bmatrix}, \quad
        V_2 = 
    \begin{bmatrix}
        -0.0270 + 0.0788\iota \\
        0.0238 + 0.0006 \iota \\
        1 
    \end{bmatrix}, \quad
        V_3 = 
    \begin{bmatrix}
        -0.0270 - 0.0788\iota \\
        0.0238 - 0.0006 \iota \\
        1
    \end{bmatrix}
\end{equation}
\noindent
Thus, the solution to linearized system is
\begin{align}
    X(t) &= c_1
    \begin{bmatrix}
    0.0022\\
    -0.0024334\\
    1
    \end{bmatrix}
    e^{(-2.6433)t} + c_2 
    \begin{bmatrix}
    -0.0270 + 0.0788\iota \\
    0.0238 + 0.0006 \iota \\
    1 
    \end{bmatrix} \\&
    e^{(-0.0191+0.0658\iota)t} +  c_3
    \begin{bmatrix}
    -0.0270-0.0788\iota\\
    0.0238-0.0006\iota\\
    1
    \end{bmatrix}
    e^{(-0.0191-0.065\iota)t}.
\end{align}
\noindent
The real valued solution corresponding state variables as defined in the linearized model are as follows
\begin{align}
    &U(t)= (0.0022)c_1e^{(-2.6433)t}+c_2e^{-0.0191)t}[-0.0270\cos(0.0658t)-0.0788\nonumber\\
    &\sin(0.0658t)]+c_3e^{(-0.0191)t}[-0.027\sin(0.0658t)+0.0788\cos(0.0658t)]\\ \nonumber\\
        &V(t)= (-0.0024)c_1e^{(-2.6433)t}+c_2e^{(-0.0191)t}[0.0238\cos(0.0658t)+0.0006\nonumber\\ &\sin(0.0658t)]+c_3e^{(-0.0191)t}[0.0238\sin(0.0658t)-0.0006\cos(0.0658t)]\\ \nonumber\\
        &W(t)= c_1e^{-2.6433t}+c_2e^{-0.0191t}\cos(0.0658t)+ c_3e^{-0.0191t}\sin(0.0658t)
\end{align}
\noindent
The solution of linearized system, as given in terms of perturbed state variables $U(t), \;V(t)$ and $W(t)$ is tested under various initial conditions and is also tested numerical in the Chapter~\ref{chapter2}.

\subsection{Conclusion}
It is concluded that when different initial conditions are given and solution of linearized system of ODEs is plotted (Figure \ref{pics:HIV_System_1_Figure_linearized}), then it is seen that curves of solution representing healthy CD4+ cells, infected CD4+ cells and virus load always reaches to the origin of the euclidean three-dimensional coordinate system, which implies that linearized system is stable around the equilibrium point and so is nonlinear system of ODEs around critical point in the HIV dynamical model without treatment $(\ref{HIV_System_11})-(\ref{HIV_System_13})$ around the equilibrium point. The drawback of linearized system is that it is only valid around the equilibrium point i.e. the solution is not considered to be the global solution.

\appendix
\section{Mathematical Background}\label{appendix-a}
This appendix present the basic mathematical concepts and definitions~\cite{Zill_ODEs_2012,Burden2010} that are used in this paper to simulate the dynamical system for HIV.
\subsection{Ordinary Differential Equation} \label{appendix-section-1}
If an equation contains ordinary derivatives of one or more unknown functions with respect to single independent variable is called ordinary differential equation. Th $n^{\mathrm{th}}$ order ordinary differential equation in the normal form is,
\begin{equation}\label{generalode}
    \frac{d^n y(t)}{dt^n} = f(t,y(t),y^{'}(t),y^{''}(t),\cdots,y^{n-1}(t))
\end{equation}
for $a \leq t \leq b$ and $t\in \Re$ where $\Re$ is set of Real Values. For example, $dx(t)/dt = 5x(t) - 3$. 
Any function $g(t)$, defined on interval $I$ possesses atleast $n$
derivatives that are continuous on $I$, which when substituted into an $n^{\mathrm{th}}$ order
ordinary differential equation reduces the equation to an identity, is said to be
a solution of the equation on the interval $I$. 
If $g(t)$ is the solution of the (\ref{generalode}) then it must satisfies the given ODE and initial conditions.
\begin{equation}
    \frac{d^n g(t)}{dt^n} - f(t,g(t),g^{'}(t),g^{''}(t),...,g^{n-1}(t)) = 0\nonumber
\end{equation}
\subsubsection{Classification of ODEs}
The ordinary differential equations (ODEs) are classified by order and linearity as follows.
\paragraph{Order} 
The order of an ordinary differential equation is the order of the
highest derivative in the equation. For example, ${d^{2}x(t)}/{dt^2} = 5x(t) - 3$
is a $2^{\mathrm{nd}}$ order ordinary differential equation.
\paragraph{Linearity} An $n^{\mathrm{th}}$ order ordinary
differential equation is said to be linear if it is of
the following form. 
\begin{equation} \label{ode_linear_equation} a_{n}(t)\frac{d^{n}x(t)}{dt^{n}}+a_{n-1}(t)\frac{d^{n-1}x(t)}{dt^{n-1}}+\cdots+a_{2}(t)\frac{d^{2}x(t)}{dt^{2}}+a_{1}(t)\frac{dx(t)}{dt} + a_{0}(t)x(t) = f(t) 
\end{equation}
where $x(t)$ has atleast $n$ derivative with respect to $t$, hence is the solution of~(\ref{ode_linear_equation}) for $t \in I$ and also $a_i(t)$ for $i=0,1,2,\cdots,n$ are continuous for every $t \in I$. 
For example, \({d^{3}x(t)}/{{dt}^3} = 5{d^2x(t)}/{dt^2} - 3\)
is $3^{\mathrm{rd}}$ order linear ordinary differential equation. 
\paragraph{Non-linear Ordinary Differential Equation}
Equations of form other than~(\ref{ode_linear_equation}) are known as Non-linear Ordinary Differential Equations. For example, $dy/dt = y^2$ for some $t$ belongs to some interval $I$ is a non-linear ODE.
\subsection{Autonomous and Non-Autonomous ODEs}\label{appendix-section-2}
An ordinary differential equation in which the
independent variable
does not appear explicitly is said to be 
\textbf{autonomous ODE}.
\begin{equation} \label{autonomous_ode}
    \frac{d^{n}y(t)}{dt^{n}} = f(y,y^{'},y^{''},y^{'''},\cdots,y^{n-1})
\end{equation}
\noindent
where $t$ belongs to some interval $I$ and $y(t)$ be the solution of~(\ref{autonomous_ode})
For example, ${dx(t)}/{dt} = x(t) + \sin(x(t))$ is a autonomous equation.  
If the any term of the an ordinary differential equation contain independent variable explicitly then it is
said to be \textbf{non-autonomous ODE.}
\begin{equation} \label{non_autonomous_ode}
    \frac{d^{n}y(t)}{dt^{n}} = f(t,y,y^{'},y^{''},y^{'''},\cdots,y^{n-1})
\end{equation}
\noindent
where $t$ belongs to some interval $I$ and $y(t)$ be the solution of~(\ref{non_autonomous_ode}), for example, ${dx(t)}/{dt} = x(t) + \sin(x(t)) + t^2 + 7t$ is a non-autonomous equation.
\subsection{Dynamical System}\label{appendix-section-3}
A dynamical system is any system that allows one to determine (at least theoretically) the future states of the system given its present or past state. For example, the recursive formula (difference equation)
$x_{n+1} = (1.05)x_n$, for $n=0,1,2,\cdots$ is a dynamical system, since the future state $x_{n+1}$ is determined given the previous state, $x_n$. 
The differential equations are also know as a dynamical system such as initial value problem (the information about the system is known at either the present or past state). For example, $dx/dt = -2x,$ where the solution $x(t)$ shows the state of the system at any time  $t$. If the initial condition is known i.e. $x(t_0) = x_0$, then system state can be determined at any past, present or future instant of time by finding the exact solution i.e. $x(t) =x_0 e^{-2(t-t_0)}$ for $t \geq t_0$.
\noindent
Following are the the types of behaviour of dynamical systems.
\begin{enumerate}
\item{\textbf{Equilibrium solutions: }A constant solution of the
autonomous differential equation \({dx}/{dt} = f(x) \) is called
an equilibrium of the equation. In other words, it is a solution
which satisfies $f(x) = 0$. The solutions either converge to the
equilibrium or diverge away from it.} 

\item{\textbf{Periodic orbits:} A periodic solution $x=x(t)$, $y=y(t)$. A periodic solution is called a
cycle. If $p$ is the period of the solution, then $X(t+p)=X(t)$ and a particle
placed on the curve at $X_0$ will cycle around the curve and return to $X_0$ in
$p$ units of time (see~\cite{Zill_ODEs_2012}, page no. $366$).}

\item{ \textbf{Chaotic orbit:} An orbit that exhibits an
unstable behavior that is not itself fixed or periodic is called
a chaotic orbit. At any point in such an orbit, there are points
arbitrarily near that will move away from the point during
further iteration. In terms of solutions, it means they are very
sensitive to small perturbations in the initial conditions and
almost all of them do not appear to be either periodic or
converge to equilibrium solutions }
\end{enumerate}
\section{MATLAB Code}\label{appendix-b}
\begin{center}
\begin{lstlisting}[frame=single]
%% main file
%% Constants
d = 0.02;
k = 100;
s = 10;
beta = 2.4e-5;
m1 = 2.4;
m2 = 0.24;
%% Initial Conditions.
ini_x = [500 800 1000 1200];
ini_y = [1e-6 10 50 0];
ini_z = [60 70 30 100];
name1 = {'T(0) = 500','T(0) = 800','T(0) = 1000','T(0) = 1200'};
name2 = {'T*(0) = 1e-6','T*(0) = 10','T*(0) = 50','T*(0) = 0'};
name3 = {'V(0) = 60','V(0) = 70','V(0) = 30','V(0) = 100'};
%% Interval of Computation
a = 0;
b = 400;
h = 0.1;  %Step size in the interval.
%% Plotting Style
CC = {'y','r','g','b','g','y',[.5 .6 .7],[.8 .2 .6]};
Markers = {'+','o','*','x','v','d','^','s','>','<'};
linespec = {'-',':','--','-.'};
%% System of Equations 2.1
f1 = @(t,x,y,z)s-d*x - beta*x*z;
f2 = @(t,x,y,z)beta*x*z - m2*y;
f3 = @(t,x,y,z)k*y - m1*z;
for i=1:length(ini_x)
    ini = [ini_x(i) ini_y(i) ini_z(i)];
    RK4_Method = RK4(a,b,h,ini,f1,f2,f3);
    h1 = figure(1)
    plot(RK4_Method(:,1),RK4_Method(:,2),linespec{i},'color',...
    CC{i},'LineWidth',2,'DisplayName',name1{i})
    axis tight
    xlabel('Time (Days)'), ylabel('T(t)')
    legend('Location','northeast')
    hold on
    %% fig file writing
    set(findall(gcf,'-property','FontSize'),'FontSize',12)
    set(gca,'FontName','Times New Roman')
    set(gcf,'Units','Inches','Position',  [0,5.8, 5.5, 2])
    h2 = figure(2)
    plot(RK4_Method(:,1),RK4_Method(:,3),linespec{i},'color',...
    CC{i},'LineWidth',2,'DisplayName',name2{i})
    axis tight
    xlabel('Time (Days)'), ylabel('T*(t)')
    legend
    hold on
    %% fig file writing
    set(findall(gcf,'-property','FontSize'),'FontSize',12)
    set(gca,'FontName','Times New Roman')
    set(gcf,'Units','Inches','Position',  [0,3, 5.5, 2])
    h3 = figure(3)
    plot(RK4_Method(:,1),RK4_Method(:,4),linespec{i},'color',...
    CC{i},'LineWidth',2,'DisplayName',name3{i})
    axis tight
    xlabel('Time (Days)'), ylabel('V(t)')
    legend('Location','northeast')
    hold on
    %% fig file writing
    set(findall(gcf,'-property','FontSize'),'FontSize',12)
    set(gca,'FontName','Times New Roman')
    set(gcf,'Units','Inches','Position',  [0,0.3, 5.5, 2])
    h4 = figure(4)
    plot3(RK4_Method(:,2),RK4_Method(:,3),RK4_Method(:,4),...
    linespec{i},'color',CC{i},'LineWidth',2)
    % title('Plotting the system')
    axis tight
    xlabel('T(t)'), ylabel('T*(t)'), zlabel('V(t)')
    grid on
    hold on
    %% fig file writing
    set(findall(gcf,'-property','FontSize'),'FontSize',12)
    set(gca,'FontName','Times New Roman')
    set(gcf,'Units','Inches','Position',  [8, 4, 5, 3])
end
%% Save the file as PNG
print(h1,'RK4_HealthyCells_System_1','-dpng','-r600');
print(h2,'RK4_InfectedCells_System_1','-dpng','-r600');
print(h3,'RK4_VirusLoad_System_1','-dpng','-r600');
print(h4,'RK4_PhasePortrait_System_1','-dpng','-r600');
\end{lstlisting}
mfile: main.m
\end{center}
\begin{center}
\begin{lstlisting}[frame=single]
function [D] = RK4(a,b,h,ini,f1,f2,f3)
% a = initial value of independent variable
% b = final value of independent variable
% h = increment in the independent varibale interval of axis.
% ini = initial conditions subject to given ode
% f1 = general function e.g. @(x,y)sin(x,y)
% f2 and f2 will also be function like f1.
N = (b-a)/h;
t = zeros(N+1,1);
x = zeros(N+1,1);
y = zeros(N+1,1);
z = zeros(N+1,1);
x(1)= ini(1);
y(1)= ini(2);
z(1)= ini(3);
t(1)= a;
for i =1:N
    t(i+1) = a + i*h;
    k1 = h*f1(t(i),x(i),y(i),z(i));
    l1 = h*f2(t(i),x(i),y(i),z(i));
    m1 = h*f3(t(i),x(i),y(i),z(i));
    k2 = h*f1(t(i)+0.5*h,x(i)+0.5*k1,y(i)+0.5*l1,z(i)+0.5*m1);
    l2 = h*f2(t(i)+0.5*h,x(i)+0.5*k1,y(i)+0.5*l1,z(i)+0.5*m1);
    m2 = h*f3(t(i)+0.5*h,x(i)+0.5*k1,y(i)+0.5*l1,z(i)+0.5*m1);
    k3 = h*f1(t(i)+0.5*h,x(i)+0.5*k2,y(i)+0.5*l2,z(i)+0.5*m2);
    l3 = h*f2(t(i)+0.5*h,x(i)+0.5*k2,y(i)+0.5*l2,z(i)+0.5*m2);
    m3 = h*f3(t(i)+0.5*h,x(i)+0.5*k2,y(i)+0.5*l2,z(i)+0.5*m2);
    k4 = h*f1(t(i)+h,x(i)+k3,y(i)+l3,z(i)+m3);
    l4 = h*f2(t(i)+h,x(i)+k3,y(i)+l3,z(i)+m3);
    m4 = h*f3(t(i)+h,x(i)+k3,y(i)+l3,z(i)+m3);
    x(i+1) = x(i) + (1/6)*(k1+2*k2+2*k3+k4);
    y(i+1) = y(i) + (1/6)*(l1+2*l2+2*l3+l4);
    z(i+1) = z(i) + (1/6)*(m1+2*m2+2*m3+m4);
end
D = [t,x,y,z];
end
\end{lstlisting}
\end{center}

\bibliographystyle{plain}
\bibliography{dissertation}

@ARTICLE{Craig2004,
author={I. K. {Craig} and {Xiaohua Xia} and J. W. {Venter}},
journal={IEEE Transactions on Education},
title={Introducing HIV/AIDS education into the electrical engineering curriculum at the University of Pretoria},
year={2004},
volume={47},
number={1},
pages={65--73},
doi={10.1109/TE.2003.817620},
ISSN={1557-9638},
month={Feb},}

@ARTICLE{Craig2005a,
  author = {I. Craig and Xiaohua Xia},
  title = {Can HIV/AIDS be controlled? Applying control engineering concepts outside traditional fields},
  journal = { IEEE Control Systems},
  year = {2005},
  volume = {25},
  pages = {80 -- 83},
  issue = {1},
  doi = {10.1109/MCS.2005.1388805},
  quality = {1}
}

@ARTICLE{Moysis2016,
  author = {Moysis, Lazaros and Kafetzis, Ioannis and Politis, Marios},
  title = {A dynamic model for {HIV} infection},
  journal = {Conference Paper},
  year = {2016},
  doi = {https://www.researchgate.net/publication/303684537},
}

@book{Nagle2011,
   title =     {Fundamentals of Differential Equations, 8th Edition},
   author =    {R. Kent Nagle, Edward B. Saff, Arthur David Snider},
   publisher = {Addison Wesley},
   isbn =      {0321747739,9780321747730},
   year =      {2011},
   series =    {},
   edition =   {8th Edition},
   volume =    {},
   url =       {http://gen.lib.rus.ec/book/index.php?md5=3b0ff6478a51cf2e0689c32d98c04660}
}

@book{Burden2010,
   title =     {Numerical Analysis, 9th Edition  },
   author =    {Richard L. Burden, J. Douglas Faires},
   publisher = {Brooks Cole},
   isbn =      {0538733519,9780538733519},
   year =      {2010},
   series =    {},
   edition =   {9},
   volume =    {},
   url =       {http://gen.lib.rus.ec/book/index.php?md5=17D673B47AA520F748534F6292F46A2B}
}

@website{khan_academy,
   title =     {{HIV} Infection Course on Khan Academy},
   publisher = {Khan Academy},
   series =    {Course Series},
   url =       {https://www.khanacademy.org/science/health-and-medicine/infectious-diseases/hiv-and-aids/a/what-is-hivaids}
}

@article{Shishana2004,
author = {Connolly, Catherine and Shisana, Olive and Shishana, Olive and Stoker, David},
year = {2004},
month = {10},
pages = {776--81},
title = {Epidemiology of {HIV} in South Africa - Results of a national, community-based survey},
volume = {94},
journal = {South African medical journal = Suid-Afrikaanse tydskrif vir geneeskunde}
}

@website{HIV_africa,
    title = {{HIV} and {AIDS} in South Africa},
    url = {https://www.avert.org/professionals/hiv-around-world/sub-saharan-africa/south-africa}
}

@website{WIKI_HIV,
    title = {Management of HIV/AIDS at Wikipedia},
    url = {https://en.wikipedia.org/wiki/Management_of_HIV/AIDS}
}

@book{Zill_ODEs_2012,
   title =     {A First Course in Differential Equations with Modeling Applications},
   author =    {Dennis G. Zill},
   publisher = {Brooks Cole},
   isbn =      {1111827052,9781111827052},
   year =      {2012},
   series =    {},
   edition =   {10},
   volume =    {},
   url =       {http://gen.lib.rus.ec/book/index.php?md5=ad896673659d5255a012c433a176c099}
}

@article{Miguel2007,
author = {Barão, Miguel and Lemos, J.M.},
year = {2007},
month = {07},
pages = {248-257},
title = {Nonlinear control of HIV-1 infection with a singular perturbation model},
volume = {2},
journal = {Biomedical Signal Processing and Control},
doi = {10.1016/j.bspc.2007.07.011}
}

@article{Mhawej2009,
author = {Mhawej, M. and Moog, Claude and Biafore, Federico},
year = {2009},
month = {01},
pages = {1263--1266},
title = {The {HIV} Dynamics is a Single Input System},
volume = {23},
journal = {IFMBE Proceedings},
doi = {10.1007/978-3-540-92841-6_310}
}

@article{Mhawej2010,
author = {Mhawej, Marie-José and Moog, Claude and Biafore, Federico and Brunet-François, Cécile},
year = {2010},
month = {01},
pages = {45--52},
title = {Control of the {HIV} infection and drug dosage},
volume = {5},
journal = {Biomedical Signal Processing and Control},
doi = {10.1016/j.bspc.2009.05.001}
}

@Article{Callaway2002,
author="Callaway, Duncan S.
and Perelson, Alan S.",
title="HIV-1 infection and low steady state viral loads",
journal="Bulletin of Mathematical Biology",
year="2002",
month="Jan",
day="01",
volume="64",
number="1",
pages="29--64",
issn="1522-9602",
doi="10.1006/bulm.2001.0266",
url="https://doi.org/10.1006/bulm.2001.0266"
}

@article{Charlotte2010,
author = {Ho, Charlotte and Ling, Bingo},
year = {2010},
month = {04},
pages = {1279--1292},
title = {Initiation of {HIV} therapy},
volume = {20},
journal = {I. J. Bifurcation and Chaos},
doi = {10.1142/S0218127410026484}
}

@article{Verica2009,
author = {Radisavljevic-Gajic, Verica},
year = {2009},
month = {04},
pages = {1251--61},
title = {Optimal Control of HIV-Virus Dynamics},
volume = {37},
journal = {Annals of biomedical engineering},
doi = {10.1007/s10439-009-9672-7}
}

@article{Kramer1999,
author = {Kramer, I.},
year = {1999},
month = {03},
pages = {95--112},
title = {Modeling the dynamical impact of {HIV} on the immune system: Viral clearance, infection, and AIDS},
volume = {29},
journal = {Mathematical and Computer Modelling - MATH COMPUT MODELLING},
doi = {10.1016/S0895-7177(99)00057-6}
}

@article{Rivadeneira2012,
author = {Rivadeneira, Pablo and Moog, Claude},
year = {2012},
month = {05},
pages = {8462--8474},
title = {Impulsive control of single-input nonlinear systems with application to {HIV} dynamics},
volume = {218},
journal = {Applied Mathematics and Computation},
doi = {10.1016/j.amc.2012.01.071}
}

@article{Rivadeneira2014,
author = {Rivadeneira, Pablo and Moog, Claude and Stan, Guy-Bart and Brunet, Cecile and Raffi, François and Ferré, Virginie and Costanza, Vicente and Mhawej, Marie-José and Ernst, Damien and Fonteneau, Raphael and Biafore, Federico and Ouattara, Djomangan and Xia, Xiaohua},
year = {2014},
month = {07},
pages = {},
title = {Mathematical Modeling of {HIV} Dynamics After Antiretroviral Therapy Initiation: A Review},
journal = {AIDS Research and Human Retroviruses},
doi = {10.1089/AID.2013.0286}
}

@article{Nowak2001,
author = {Nowak, Martin and May, Robert},
year = {2001},
month = {01},
pages = {},
title = {Virus Dynamics: Mathematical Principles of Immunology And Virology}
}

@article{Landi2008,
author = {Landi, Alberto and Mazzoldi, Alberto and Andreoni, Chiara and Bianchi, Matteo and Cavallini, Andrea and Laurino, Marco and Ricotti, Leonardo and Iuliano, Rodolfo and Matteoli, Barbara and Ceccherini-Nelli, Luca},
year = {2008},
month = {03},
pages = {162--8},
title = {Modelling and control of {HIV} dynamics},
volume = {89},
journal = {Computer methods and programs in biomedicine},
doi = {10.1016/j.cmpb.2007.08.003}
}

@article{Xiaohua2006,
author = {Xia, Xiaohua},
year = {2006},
month = {12},
pages = {485--492},
title = {MODELLING OF {HIV} INFECTION: VACCINE READINESS, DRUG EFFECTIVENESS AND THERAPEUTICAL FAILURES},
volume = {39},
journal = {IFAC Proceedings Volumes},
doi = {10.3182/20060402-4-BR-2902.00485}
}

\end{document}